\RequirePackage{fix-cm}
\documentclass[smallextended]{svjour3}       
\smartqed  
\usepackage{amsmath}
\usepackage{latexsym}
\usepackage{hyperref}
\usepackage{amssymb}
\usepackage{lipsum}
\usepackage{amsfonts}
\usepackage{graphicx}
\usepackage{epstopdf}
\usepackage{amsopn}
\usepackage{color}
\usepackage{xcolor}

\newcommand{\esp}{\mathbb{E}}
\newcommand{\espn}{\mathbf{E}}
\newcommand{\prob}{\mathbb{P}}
\newcommand{\probn}{\mathbf{P}}
\newcommand{\var}{\mathbb{V}}

\newcommand{\re}{\mathbb{R}}

\newcommand{\ignore}[1]{}
\newcommand{\cbdg}{{\bf c}_{\sf bdg}}

\DeclareMathOperator*{\gap}{gap}

\newcommand{\vertiii}[1]{{\left\Vert\kern-0.4ex #1 
\kern-0.4ex\right\Vert}}

\newcommand{\Lpnorm}[1]{\vertiii{\,#1\,}_{p}}

\newcommand{\Lunorm}[2]{\vertiii{\,#1\,}_{#2}}

\renewcommand{\c}[1]{\mathcal{#1}}

\newcommand{\dist}{\mathsf{d}}
\newcommand{\diam}{\mathsf{diam}}

\newcommand{\cdiam}{3\sqrt{2}}
\newcommand{\cgamma}{2\sqrt{2}}
\newcommand{\cgammaconc}{4\sqrt{3}}
\newcommand{\cdiamconc}{6\sqrt{3}}
\newcommand{\stheta}{\vartheta}

\newcommand{\parametersgoodness}{(\sigma_*^2,\sigma^2,\alpha,\rho)}
\newcommand{\parametersgreatness}{(\sigma_*^2,\sigma^2,\alpha,p,\kappa_p)}

\newcommand{\parametersgoodnessplus}{(\sigma_*^2,\sigma^2(\stheta,\delta),\alpha,\rho)}
\newcommand{\parametersgreatnessplus}{(\sigma_*^2,\sigma^2(\stheta,\delta),\alpha,p,\kappa_p)}

\newtheorem{assumption}{Assumption}

\begin{document}

\title{Sample average approximation with heavier tails I
}
\subtitle{Non-asymptotic bounds with weak assumptions and stochastic constraints}


\author{Roberto I. Oliveira         \and
        Philip Thompson 
}


\institute{Roberto I. Oliveira \at
              Instituto de Matem\'atica Pura e Aplicada (IMPA), Rio de Janeiro, RJ, Brazil. \\
              \email{rimfo@impa.br}           
           \and
           Philip Thompson \at
           Purdue University \& Krannert School of Management, West Lafayette, USA. \\
           \email{thompsp@purdue.edu}
}

\date{Received: date / Accepted: date}

\maketitle

\begin{abstract}
We derive new and improved non-asymptotic deviation inequalities for the sample average approximation (SAA) of an optimization problem. Our results give strong error probability bounds that are ``sub-Gaussian"~even when the randomness of the problem is fairly heavy tailed. Additionally, we obtain good (often optimal) dependence on the sample size and geometrical parameters of the problem. Finally, we allow for random constraints on the SAA and unbounded feasible sets, which also do not seem to have been considered before in the non-asymptotic literature. Our proofs combine different ideas of potential independent interest: an adaptation of Talagrand's ``generic chaining"~bound for sub-Gaussian processes; ``localization"~ideas from the Statistical Learning literature; and the use of standard conditions in Optimization (metric regularity, Slater-type conditions) to control fluctuations of the feasible set. 
\end{abstract}

\section{Introduction}

Understanding {\em sample average approximations} is a fundamental problem in Stochastic Programming \cite{roemisch2003,shapiro:dent:rus}.  Suppose we are given an optimization problem:
\begin{eqnarray}
f^*:=\min_{x\in Y} &&\quad f_0(x)\nonumber\\
\mbox{s.t.}&&\quad f_i(x)\le0,\quad\quad\forall i\in\mathcal{I},\label{problem:min}
\end{eqnarray}
with $Y\subset\re^d$ and a nonempty \emph{feasible set} 
\begin{eqnarray}
X:=\left\{x\in Y:f_i(x)\le0, \forall i\in\mathcal{I}\right\}.
\label{equation:feasible:set}
\end{eqnarray}

In this paper, each of the functions $f_i$ is given by an expectation
\begin{align}
f_i(x) = \espn\,F_i(x,\cdot) := \int_{\Xi}\,F_i(x,\xi)\,\probn(d\xi)\label{eq:mean:fi}
\end{align}
where $\probn$ is a probability measure over a set $\Xi$ and the $F_i:Y\times \Xi\to \re$ are Carath\'{e}odory functions. In typical settings, the measure $\probn$ and the functions $f_i$ are not directly accessible. It may be, however, that a {\em random sample} $\{\xi_k\}_{k=1}^N$ from $\probn$ is available. If that is the case, it is natural to consider the sample-average approximation (or SAA) to (\ref{problem:min}), where the $f_i$ are replaced by sample averages:
\begin{align}
\widehat{F}_i(x):= \widehat{\espn}\,F_i(x,\cdot) = \frac{1}{N}\sum_{k=1}^NF_i(x,\xi_k)\label{eq:mean:hat:fi}
\end{align}
This leads to some natural questions  considered in numerous works in stochastic optimization:
\begin{enumerate}
\item Are (nearly) optimal solutions to the SAA also nearly feasible and nearly optimal for the original problem (\ref{problem:min})?
\item Are the values of the two problems typically close?
\end{enumerate}

\emph{Asymptotic} analyses of the SAA assume the sample size $N$ diverges whereas the functions $F_i$, $f_i$, the set $Y$ and the measure $\probn$ remain fixed. Using tools such as {uniform} Strong Law of Large Numbers and Central Limit Theorems, these analyses obtain precise answers to the above questions. This program has been carried out in numerous works, e.g.,  \cite{artstein:wets,dupacova:wets,king:rockafellar,king:wets,pflug1995,pflug1999,pflug2003,shapiro1989,shapiro1991,shapiro2003,shapiro:dent:rus}. See \cite{shapiro:dent:rus,tito:bayraksan,pasupathy} for extensive reviews.

Another type of analysis, which we pursue in this work, is \emph{non-asymptotic} in nature. It consists of proving explicit bounds for the value and quality of SAA solutions with explicit dependence on the sample size $N$ and other problem parameters. For instance, letting $f^*$ and $\widehat{F}^*$ be the values of the original problem (\ref{problem:min}) and its SAA (respectively), a recent non-asymptotic result by Guigues, Juditsky and Nemirovski \cite{guigues:juditsky:nemirovski} gives guarantees of the form:
 \begin{equation}\label{eq:guarantees}\forall t\geq 0\,:\,\prob\left\{|\widehat{F}^* - f^*|\leq \frac{A + B\sqrt{t}}{\sqrt{N}}\right\} \geq 1 - e^{-t},\end{equation}
 where $A$ and $B$ do not depend on $N$ or $t$ (but do depend on other problem parameters).  Guarantees of this kind are called ``sub-Gaussian"\footnote{Another typical light-tail condition is to assume an sub-exponential tail.} because they imply that the tail decay $\sqrt{N}|\widehat{F}^* - f^*|$ roughly matches that of a Gaussian distribution with standard deviation $B$. This sort of asymptotic behavior is what one expects from asymptotic statements such as those found in \cite{shapiro1989}.

With few exceptions, non-asymptotic guarantees in the literature require that the random variables $F_i(x,\xi)$ be very light-tailed: that is $|F_i(x,\xi)-f_i(x)|$ has finite $p$-th moments for all $p\geq 1$. In the rare cases where this is avoided \cite{kancova:omelchenko2015}, the dependence on $N$ is suboptimal (as we shall see). Other limitations to current finite-sample analyses of SAA include requiring the feasible set $X$ to be bounded, and avoiding expected value constraints. Even in an equation like (\ref{eq:guarantees}), it is often not clear if the dependence of ``constants"~like $A$ and $B$ on other problem parameters (such as the dimension) is reasonable. 

There is thus a gap between what one may expect SAA to do on the basis of asymptotic analyses, and what has been proven to do non-asymptotically. Is this a technical issue, or does it point to underlying limitations of SAA? This question is especially pressing in high-dimensional problems, where asymptotic theory is not expected to give good results even for fairly large $N$. Luckily, one {\em can} prove significantly better finite-sample guarantees for SAA, as we explain below.

\subsection{Our contribution}

Our goal in this paper is to obtain new and improved non-asymptotic bounds for the sample average approximation. Our probabilistic assumptions are significantly weaker than in previous work, and our bounds often improve on other results by making better use of the geometry of our problem. We highlight some salient features of our approach. 

\paragraph{Finite-moment assumptions.} We {\em do not require} infinitely many moments of any of the random variables involved in our problem. Our main assumption is that, given norm $\Vert\cdot\Vert$ over $\re^d$, the $F_i:Y\times \Xi\to\re$ are stochastically H\"{o}lder over sets $Z\subset Y$, in the sense that inequalities of the kind
\begin{equation}\label{eq:continuity}\forall x,x'\in Z\,:\,|F_i(x,\xi) - F_i(x',\xi)|\leq \mathsf{L}_i(\xi)\|x-x'\|^\alpha
\end{equation}
hold in suitable $Z\subset Y$, with $0<\alpha\leq 1$ and $\mathsf{L}_i(\xi)$ satisfying weak conditions. See Assumptions \ref{assump:goodFi} and \ref{assump:greatFi} below for details. Conditions of the kind of (\ref{eq:continuity}) have often appeared in the literature \cite{kancova:houda2015,kancova:omelchenko2015}, but either with much stronger moment assumptions on the $\mathsf{L}_i$ or with suboptimal error bounds in the sample size.

\paragraph{Joint guarantees for values, feasibility and optimality.} An inequality such as (\ref{eq:guarantees}) bounds the difference in values between the SAA and the original problem. Our results also quantify how good the extent to which SAA is close to being feasible and optimal for the original problem. In all cases, we obtain optimal dependence on the same size $N$, as well as ``sub-Gaussianity"~for a relevant set of parameters. In Section \ref{sub:discussion}, we comment on what parameters we consider. Their precise definitions are discussed in more detail in Sections \ref{sec:results} and \ref{sec:convexresults} (see in particular the discussion in Section \ref{ss:prelim:discussion}). In Section \ref{sec:results}, Theorem \ref{thm:generalrandomset} considers general (possibly non-convex) problems. In Section \ref{sec:convexresults},  Theorem \ref{thm:convexrandomset} and Propositions \ref{prop:LRSS}-\ref{prop:LRSS:random:set} state sharper ``localized'' bounds for convex optimization.

\paragraph{Generic chaining without light tails.} A key step in our proofs will be to obtain concentration inequalities for $\widehat{\espn}F_i(x,\cdot) - f_i(x)$ under assumptions such as (\ref{eq:continuity}). For this purpose, we adapt to our setting Talagrand's generic chaining method for empirical processes \cite{talagrand1994,talagrand2014}, as improved by Dirksen \cite{dirksen2015}. Generic chaining is an optimal method for taking problem geometry into account, and gives good problem-dependent bounds on ``constants"~like $A$ and $B$ in (\ref{eq:guarantees}) under sub-Gaussian assumptions. We obtain novel concentration generic chaining inequalities (Theorem \ref{thm:concentration}) that do not require light tails, which are of independent interest. 

\paragraph{Localization, convexity and unbounded sets.} ``Localization"~is a key idea developed by researchers in Statistical Learning, especially Koltchinskii, Mendelson and their collaborators \cite{koltchinskii2000,koltchinskii2006,bartlett2005,%
bartlett2006,mendelson2015,mendelson2017}. For convex problems, it means that ``failure"~for an SAA solution must originate from ``bad behavior"~of the SAA in a (often small) sublevel set around the minimum. We will show that this idea often leads to faster convergence rates for SAA. It also allows us to only require the H\"{o}lder condition (\ref{eq:continuity})  in a potentially ``small"~subset $Z\subset Y$. In some cases, this allows us to consider unbounded convex feasible sets and functions $F_i$ with superlinear growth. Theorem \ref{thm:convexrandomset}
presents general localized rates.  Propositions \ref{prop:LRSS} and \ref{prop:LRSS:random:set} exemplify Theorem \ref{thm:convexrandomset} when typical regularity assumptions hold.

\paragraph{Constraints in expectation.} We deal systematically with constraints in expectation. These mean that the feasible set of the SAA is a perturbation of a deterministic set. We control these perturbations by combining tools from Optimization theory -- metric regularity and Slater-type conditions -- with our ``localization toolbox".  A key result will be to show that, when constraints are perturbed, this does not change much the ``generic chaining"~parameters of relevant sub-level sets of the objective function. We remark that we do not make detailed reference to the large literature on optimization with \emph{chance constraints}. This challenging problem is out of scope of this paper as the continuity assumption in \eqref{eq:continuity} is not satisfied.

\paragraph{Examples.} Finally, we present four different applications of our theory. The first two examples is treated in detail in this work as a proof of concept of Theorems \ref{thm:generalrandomset} and \ref{thm:convexrandomset}.  The other examples, which require finer analyses, are presented in a dedicated companion paper \cite{2020oliveira:thompsonII}.

\begin{example}[Regular convex optimization problems; Section \ref{ss:particular:case}]
We consider SAA with convex objective and constraints satisfying two typical regularity conditions: (1) a local Slater constraint qualification (Assumption \ref{assump:LSCQ}) and (2) a local regular solution set (Assumption \ref{assump:LRSS}). The first is typical while 
the latter is satisfied, e.g., for objectives that are locally strongly convex or with local weak sharp minima \cite{2005burke:deng,1993burke:ferris}. We consider constraint-free problems (Proposition \ref{prop:LRSS}) or problems with random constraints (Proposition \ref{prop:LRSS:random:set}). They offer concrete localized rates implied by the general Theorem \ref{thm:convexrandomset}. In particular, unlike Theorem \ref{thm:generalrandomset}, the obtained rates depend only on the diameter and a complexity measure of a neighbourhood of the solution set. See futher discussions in Sections \ref{ss:prelim:discussion} and \ref{ss:particular:case}.  
\end{example}

\begin{example}[Metric projection problems; Section \ref{sec:simpleresults}] We consider the special case of problem \ref{problem:min} where the feasible set is convex and $f_0(x):=\|x-x_0\|_2^2$ with $x_0\in\re^d$ fixed and $\|\cdot\|$ the standard Euclidean norm. In this case, the unique optimal solution of our problem is the metric projection of $x_0$ onto the feasible set $X$.  We provide finite-sample guarantees for SAA that make strong use of localization. One particular difficulty of this problem is the fact that the Lipschitz modulus of the objective varies along the feasible set.
\end{example}

\begin{example}[Risk-averse portfolio optimization; in companion paper \cite{2020oliveira:thompsonII}]\label{ex:CVaR} Here, \[\xi = (\xi[1],\dots,\xi[d])^T\] is a random vector whose coordinates correspond to losses of $d$ distinct financial assets. If $x=(x[1],\dots,x[d])^T$ is a vector whose coordinates describe the fractions of the initial capital invested in assets $1,\ldots,d$, then the total loss is proportional to $\langle x,\xi\rangle$. We wish to minimize the expectation of $\langle x,\xi\rangle$ subject to a constraint on the conditional value-at-risk of the solution \cite{rockafellar:urysaev2000}. In this problem, the case of light-tailed $\xi$ would be of little interest. In a companion paper \cite{2020oliveira:thompsonII}, we describe specific assumptions that allow for heavy tails. 
We show that the localization toolbox obtained in this paper implies that ``risk inflation'' only affects a lower dimensional space.
\end{example}

\begin{example}[The Lasso estimator; in companion paper \cite{2020oliveira:thompsonII}]
In Least-Squares-type problems, the loss function to be minimized is $f(x)=\espn\,F(x,\cdot)$, with
$
F(x,\xi):=\left[y(\xi)-\langle \mathbf{x}(\xi),x\rangle\right]^2.
$
Here, $y(\xi)\in\re$ and the random vector $\mathbf{x}(\xi)\in \re^d$. Minimizing the empirical function $\widehat{F}(x):=\widehat{\espn}\,F(x,\cdot)$ tend to work when $N\gg d$, but not when $N\ll d$, as the problem is undetermined. Tibshirani \cite{tibshirani1996} proposed the Lasso estimator given by the problem
$
\min_{x\in Y}\widehat F(x),
$
with
$
Y:=\{x\in\re^d:\Vert x\Vert_1\le R\},
$
where $R>0$ is a tuning parameter and $\Vert\cdot\Vert_1$ denotes the $\ell_1$-norm. Inspired by Bickel, Ritov and Tsybakov \cite{bickel:ritov:tsybakov2009}, we analyse in our companion paper \cite{2020oliveira:thompsonII} the least squares problem subjected to $\Vert\mathbf{\widehat D}_{2}x\Vert_1\le R$, where 
$\mathbf{\widehat D}_{2}$ is a data-driven matrix. We obtain improved  ``persistence'' bounds \cite{bartlett:mendelson:neeman2012}  for a least-squares Lasso-type estimator. Our proof is based on localization techniques established in this paper.
\end{example}

\subsection{Discussion and comparison with previous work} \label{sub:discussion}

Of the numerous papers on the topic of SAA, we highlight \cite{pflug1999,roemisch2003,pflug2003,atlason:epelman:henderson2004,kleigweit:shapiro:tito2001,shapiro:tito2000,%
shapiro:nemirovski2005,wang:ahmed2008,shapiro2003,shapiro:dent:rus,vogel1988,vogel2008,vogel2008:2,kancova:omelchenko2015,kancova:houda2015,guigues:juditsky:nemirovski,banholzer:fliege:werner} as relevant to our findings. Except for \cite{royset2012,kancova:houda2015,kancova:omelchenko2015}, all of the non-asymptotic papers papers assume light-tailed data. This restriction is lifted in references \cite{kancova:houda2015,kancova:omelchenko2015}, but at the cost of worse dependence on $N$ in (\ref{eq:guarantees}): the error $\widehat{F}^*-f^*$ is stochastically bounded by a quantity that decays like $N^{-\beta}$ for some $\beta<1/2$. By contrast, the paper \cite{royset2012} makes weak probabilistic assumptions on the data, and obtains distributional results in a asymptotic setting for a {\em reformulation} of the SAA. Our results assume heavy-tails, are nonasymptotic, do not use reformulations and achieve the optimal rate $N^{-\frac{1}{2}}$ in terms of the sample size, with joint guarantees of feasibility and optimality. In addition, our bounds explicitly account for the geometry of the feasible set.

To understand our improvements, it is necessary to take a step back and understand {\em how} light tails were used in previous analysis. For the moment, consider the case where $\mathcal{I}=\emptyset$, i.e. the feasible set of our original problem (\ref{problem:min}) is $X=Y$ and there are no constraints in expectation. The easiest way to bound the difference between $\widehat{F}^*-f^*$ (say) is via a uniform bound:
\begin{equation}\label{eq:uniformtailbound}|\widehat{F}^* - f^*|\leq \sup_{x\in Y}|\widehat{F}(x) - f(x)| = \sup_{x\in X}\left|\frac{1}{N}\sum_{k=1}^N(F(x,\xi_k) - \espn F(x,\cdot))\right|.\end{equation}
To bound the right hand side (RHS), a typical approach uses two steps. The first one is to discretize the feasible set $X$; this reduces the problem of controlling the supremum over $X$ to controlling the supremum over finite subsets. The next step is to use concentration-of-measure inequalities \cite{boucheron:lugosi:massart2013} to deal with the finite subsets. For this it is essential to have strong concentration bounds, which typically require light tails. Our approach uses ideas that seem new in this setting. We discretize~via Talagrand's generic chaining method \cite{talagrand2014}, which is optimal in Gaussian processes and gives better dependence on the geometry of the problem. We do this via a novel concentration inequality (Theorem \ref{thm:concentration}) that separates the fluctuations of the RHS into two components: one that is {\em always sub-Gaussian}, and another that depends on the fluctuations of $\mathsf{L}^2(\cdot)$. This will give us sub-Gaussian results in certain probability regimes. 

To continue with our approach, we note that the bound (\ref{eq:uniformtailbound}) is often too pessimistic. Oftentimes, one can show that the minimizer of the SAA is usually quite close to the minimizer of the original problem. If that is the case, then $|\widehat{F}^* - f^*|\leq \sup_{x\in Z}|\widehat{F}(x) - f(x)|$ for a potentially much smaller set $Z\subset Y$. This localization idea goes back at least to the work of Koltchinskii and Panchenko \cite{koltchinskii2000} and was more fully developed in Koltchinskii's IMS Medallion Lecture \cite{koltchinskii2006}. Mendelson has also greatly contributed to this approach, starting with joint work with Bousquet and Bartlett \cite{bartlett2005} and continuing with his papers \cite{mendelson2015,mendelson2017}. These works employed localization in a somewhat different form from \cite{koltchinskii2000,koltchinskii2006} in convex settings. In this paper, we apply and extend those ideas to the setting where there are constraints in expectation. In Proposition \ref{lem:convexlocalized}, we give a ``localized bound" for perturbations of the original problem. These include perturbations of the constraints. We then show in Lemma \ref{lem:robinson2} that one can control the effect of those perturbations on the feasible sets via Slater-type conditions.

Before concluding this section, it is instructive to discuss beforehand what are the relevant parameters appearing in the improved rates of  Theorems \ref{thm:generalrandomset}-\ref{thm:convexrandomset} and Propositions \ref{prop:LRSS}-\ref{prop:LRSS:random:set}. Let $c>0$ denote an absolute constant. For light-tailed H\"older functions with exponent $\alpha$ and modulus $\sigma$, a reanalysis of the arguments in \cite{guigues:juditsky:nemirovski} shows that the parameters in bound \eqref{eq:guarantees} are typically of the form \[A_{\eqref{eq:guarantees}}=c\sigma\diam^{\alpha}(Y)\sqrt{d}\mbox{ and }B_{\eqref{eq:guarantees}}=c\sigma\diam^{\alpha}(Y)\] where $d$ is the dimension and $\diam(Y)$ denotes the diameter of $Y$.

The main improvement of Theorem \ref{thm:generalrandomset} is to allow heavier tails and give joint guarantees for optimality and feasibility, with  bounds that are of the form \[A=c\sigma\gamma^{(\alpha)}(Y)\mbox{ and }B=c\sigma\diam^\alpha(Y)+c\sigma_*\] where $\sigma_*$ is the variance at a solution. Here $\gamma^{(\alpha)}(Z)$ denotes a complexity measure of a set $Z\subset\re^d$ coming from the theory of Gaussian processes, which we discuss in \ref{ss:complexity:parameter}. A conservative upper bound  $c\diam^{\alpha}(Z)\sqrt{d}$ is possible. In case of random constraints, the probability bounds depend logarithmically on the number of constraints and (implicitly) on the metric regularity constant $\mathfrak{c}$ of the feasible set (Assumption \ref{assump:MRF}).

Theorem \ref{thm:convexrandomset}  and Propositions \ref{prop:LRSS}-\ref{prop:LRSS:random:set} give sharper localized bounds for convex problems satisfying a Slater condition (Assumption \ref{assump:LSCQ}). The statement of Theorem \ref{thm:convexrandomset} is more involved. Qualitatively, the rates depend on factors of the form \[A=c\sigma(\epsilon,\delta)\gamma^{(\alpha)}(X^{*,\epsilon}_0)\mbox{ and } 
B=c\sigma(\epsilon,\delta)\diam^{\alpha}(X^{*,\epsilon}_0)+c\sigma_*\] where $\sigma(\epsilon,\delta)$ denotes the Holder modulus variance over the set $X^{*,\epsilon}_\delta$ of approximate solutions having feasibility slackness $\delta>0$ and optimality slackness $\epsilon>0$. Hence we allow the H\"older modulus to vary across bounded regions. ``Localization'' stems from the fact that the diameter and complexity of $X^{*,\epsilon}_\delta$ are typically much smaller than the ones for $X$ or $Y$. In case of random constraints, the range of $(\epsilon,\delta)$ for which this bound holds depend on the parameters of the Slater condition (Assumption \ref{assump:LSCQ}). We refer to Section \ref{ss:prelim:discussion} for a qualitative discussion on these points before Theorem \ref{thm:convexrandomset} is presented formally. Technical rate statements also appear in the literature on localization in Statistics and Machine Learning \cite{koltchinskii2006,mendelson2017} (in this setting without random constraints). The difficulty lies in the fact that a precise rate depends on ``solving'' a fixed-point on $(\epsilon,\delta)$. In our case, an additional difficulty is that  $(\epsilon,\delta)$ are coupled: feasibility affects optimality.

Propositions \ref{prop:LRSS}-\ref{prop:LRSS:random:set} presents specific rates implied by the general Theorem \ref{thm:convexrandomset} assuming, besides Assumption \ref{assump:LSCQ}, a typical local regularity assumption on the solution set (Assumption \ref{assump:LRSS}). This includes, e.g., cases when the objective is locally strongly convex or it has locally weakly sharp minima \cite{1993burke:ferris,2005burke:deng}.  For simplicity, we assume the H\"older modulus's variance $\sigma$ is constant. For the sake of comparison with the literature on localization in Statistical Learning, Assumption \ref{assump:LRSS} is an analog (with proper differences) of the so called local \emph{Bernstein condition} on the loss function \cite{mendelson2017}. For strong regular sets ($\kappa=1/2$),  Proposition \ref{prop:LRSS} presents ``fast-$1/N$-rates'' for the constraint-free case of the form $\mathfrak{c}\sigma^2(\mathfrak{C}_\alpha+t)/N$ where $\mathfrak{c}$ is a condition number (Assumption \ref{assump:LRSS}). Here, $\mathfrak{C}_\alpha$ is the ratio comparing the complexity and diameter of the approximate solution set. A pessimistic bound for $\mathfrak{C}_\alpha$ is of order $d$. Proposition \ref{prop:LRSS:random:set}, allowing random constraints, presents ``slower-$1/\sqrt{N}$-rates'' of the form $\mathfrak{C}\sqrt{\mathfrak{C}_\alpha+\log m+t)/N}$. Here, $m$ is the number of constraints and $\mathfrak{C}$ is a constant depending polynomially on the regularity constants of Assumptions \ref{assump:LSCQ} and \ref{assump:LRSS}, $\alpha$, $(\sigma,\sigma_*)$ and the diameter of the approximate solution set $X^{*,\vartheta_*}_0$ for a small slack $\vartheta_*>0$. A notable fact is that these rates are localized in that they do not depend on the size and complexity of the entire feasible set $X$, just of an approximate solution set. Another notable fact is the deterioration of order $\sqrt{N}$ in the rate when random constraints are present. As explained for the metric projection problem, this feature is in general unavoidable.

Finally, we emphasize a few points about our approach. In most cases, we expect our results to be of optimal or nearly optimal {\em order of magnitude} in terms of problem geometry and/or sample size. Like with most non-asymptotic analyses, we do {\em not} expect our results to be as tight as asymptotic results when it comes to constants. The goal of our paper is {\em not} to give bounds that can be directly used in practice, but rather to better understand the fundamental properties of SAA with finite samples, in settings where other problem parameters (such as the dimension and diameter) can be large.

\subsection{Organization}

The remainder of the paper is organized as follows. Section \ref{section:notation} fixes notation and recalls some notions from Probability theory, most notably ``generic chaining". Section \ref{sec:basic} presents the setup for our problem and the assumptions we require on the random variables and feasible sets involved. Section \ref{sec:results} contains the statement of our main results for possibly non-convex problems (Theorem \ref{thm:generalrandomset}). Section \ref{sec:convexresults} states our main results for convex problems (Theorem \ref{thm:convexrandomset} and Propositions \ref{prop:LRSS}-\ref{prop:LRSS:random:set}). These are immediately applied to a simple example in Section \ref{sec:simpleresults}.

The next three sections presents our main technical tools separately, as we believe they might be applied or combined in different ways. The concentration inequality for heavy-tailed distributions is presented in Section \ref{sec:concentration}. Section \ref{sec:deviationlocalization} describes the relationship between good approximation properties of the SAA and differences $\widehat{F}_i(x)-f_i(x)$. In particular, this is where we prove our localization results. It will be clear that we need to understand how the feasible set $X$ changes when constraints are slightly relaxed. We present our geometrical tools for that purpose in Section \ref{sec:perturbation}. Sections \ref{sec:deviationlocalization}-\ref{sec:perturbation} are deterministic results and may be useful elsewhere.

The paper ends with Section \ref{sec:proofs}, where the main results are proven. An Appendix presents a few technical proofs left over from the main text.  

\section{Preliminaries}\label{section:notation}

\subsection{Basic notation}  

Given a set $S$, we denote its (potentially infinite) cardinality by $|S|$. The complement of an event $E$ in a probability space is $E^c$. For $m\in\mathbb{N}$, we write $[m]:=\{1,\ldots,m\}$. 

Elements of $\re^d$ are column vectors. Given $x\in\re^d$, its coordinates are denoted by $x[i]$, $1\leq i\leq d$. A superscript $T$ is used to denote transposition of a vector, so $x\in\re^d$ is given by $(x[1],\dots,x[d])^T$. The inner product of $x,y\in\re^d$ is denoted by $\langle x,y\rangle$ or $x^Ty$. Norms are denoted by $\|\cdot\|$ and the unit ball around $0$ in that norm is $\mathbb{B}$. Given $a\in\re$, $a_+:=\max\{a,0\}$.

Let $(\mathcal{M},\dist)$ be a metric space. We let $\diam(A)$ denote the (potentially infinite) diameter of $A\subset \mathcal{M}$. Given $x\in \mathcal{M}$ and $A\subset \mathcal{M}$ nonempty, $\dist(x,A):=\inf_{a\in \mathcal{A}}\dist(x,a)$. 

We fix from now on a probability space $(\Omega,\c{A},\prob)$ and assume all random variables we consider are defined on it. Given a random variable $Z$, we let $\esp[Z]$ denote its mean, $\var[Z]$ denote its variance and $\Lpnorm{Z}:= (\esp[|Z|^p])^{1/p}$ denotes $L^p$ norm (for $p\geq 1$). 

\subsection{Complexity parameters for sets} \label{ss:complexity:parameter}

We review in this section some definitions and results about ``generic chaining". Talagrand's book \cite{talagrand2014} is the best reference for these concepts.

The ``generic chaining"~functional of a metric space $(\mathcal{M},\dist)$ is a measure of the ``complexity"~of discretizing $\mathcal{M}$ at different scales. To define it, we need the following concept. A sequence $\{\mathcal{A}_j\}_{j=0}^{+\infty}$ is {\em admissible} if each $\mathcal{A}_j$ is a partition of $\mathcal{M}$, with $|\mathcal{A}_0|=1$ and $|\mathcal{A}_j|\leq 2^{2^j}$ for each $j\geq 1$. For each $j$, we let $\diam(\mathcal{A}_j)$ to denote the largest diameter of a set in partition $\mathcal{A}_j$. 

Given $0<\alpha\leq 1$, $\gamma_{2}^{(\alpha)}(\mathcal{M},\dist)$ is defined as:
\begin{equation}\label{eq:defgamma2}\gamma_{2}^{(\alpha)}(\mathcal{M},\dist):=\inf\limits_{\{\mathcal{A}_j\}_{j}\text{ admissible}}\left\{\,\sum_{j\geq 0}2^{\frac{j}{2}}\diam(\mathcal{A}_j)^{\alpha}\right\}.\end{equation}

\begin{remark}In the usual definition of the $\gamma_2$ functional, one takes $\alpha=1$. $\gamma_{2}^{(\alpha)}$ is the functional obtained when the metric $\dist$ is replaced by the equivalent metric $\dist^{\alpha}$. We will omit $\dist$ from the notation when it is clear which metric we are referring to. This remark should be kept in mind when reading Theorem \ref{theorem:gc} and equation (\ref{eq:entropybound}) below.\end{remark}

 Talagrand's celebrated majorizing measures theorem \cite{talagrand1994,talagrand2014} shows that:
\begin{equation}\label{eq:talagrand}c\,\esp\left[\sup_{x\in \mathcal{M}}|Y_x-Y_{x_0}|\right]\leq \gamma_{2}^{(\alpha)}(\mathcal{M},\dist)\leq C\,\esp\left[\sup_{x\in \mathcal{M}}|Y_x-Y_{x_0}|\right],\end{equation}
with $c,C>0$ universal, when the $Y_x$ are mean-zero Gaussian and $\var[Y_x-Y_{x'}] =\dist(x,x')^{2\alpha}$.

In fact, the upper bound in (\ref{eq:talagrand}) does not require that the $Y_x$ be truly Gaussian, only that they have sub-Gaussian tails. The next theorem, which we will use later, illustrates this point. It follows from Talagrand's work \cite{talagrand1994,talagrand2014} with an improvement due to Dirksen \cite{dirksen2015}\footnote{The constants appearing in our Theorem \ref{theorem:gc} are not the same as in \cite{dirksen2015}, but can be easily obtained via the same method.}.

\begin{theorem}[Generic chaining tail bound \cite{talagrand1994,talagrand2014,dirksen2015}]\label{theorem:gc} Suppose $\gamma^{(\alpha)}_2(\mathcal{M})<+\infty$ (in particular, $\mathcal{M}$ is totally bounded). Let $\{Y_x\}_{x\in \mathcal{M}}$ be a family of random variables indexed by the points of $\mathcal{M}$, which depend almost surely continuously on $x$. Assume further that the $Y_x$ satisfy the following sub-Gaussian assumption.
\[\forall t\geq 0\,\forall x,x'\in \mathcal{M}\,:\, \prob\{Y_x - Y_{x'}\geq \dist(x,x')^\alpha\sqrt{2(1+t)}\}\leq e^{-t}.\]
Then for any $t\geq 0$ and $x_0\in \mathcal{M}$
\[\prob\left\{\sup_{x\in \mathcal{M}}\,|Y_x - Y_{x_0}| \leq \cdiam \diam(\mathcal{M})^\alpha\sqrt{1+t} + \cgamma\,\gamma_{2}^{(\alpha)}(\mathcal{M})\right\}\leq e^{-t}.\]\end{theorem}

The functional $\gamma^{(\alpha)}_2(\mathcal{M},\dist)$ is somewhat mysterious, and can be quite difficult to compute.  In the case 
$\mathcal{M}\subset \mathbb{R}^d$, $\alpha=1$ and $\dist$ is given by the standard Euclidean norm, Talagrand's general theory connects $\gamma^{(1)}_2(
\mathcal{M})$ to a parameter called the {\em Gaussain width}. Letting $g\in \mathbb{R}^d$ denote a standard Gaussian random vector, the Gaussian width of $\mathcal{M}$ is defined as 
\[w(\mathcal{M}):=\mathbb{E}\sup_{x\in\mathcal{M}}\langle x,g\rangle.\]
It follows from \eqref{eq:talagrand} that the ratio $\gamma_2^{(1)}(\mathcal{M})/w(\mathcal{M})$ is upper and lower bounded by absolute constants $c,C>0$. One consequence of this fact is that, if $\mathcal{M}$ is the convex hull of a finite set $F$ of points, then:
\begin{equation}\label{eq:gamma2convexhull}\gamma_2^{(1)}(\mathcal{M})\leq C'\,\sqrt{\log|F|}\,\max_{x\in F}\|x\|,\end{equation}
for an absolute constant $C'>0$.

A more general upper bound for $\gamma^{(\alpha)}_2(\mathcal{M})$ comes from Dudley's {\em entropy integral} \cite{talagrand2014}. Recall that an $r$-net in $\mathcal{M}$ is a set $A\subset \mathcal{M}$ such that $\dist(x,A)\leq r$ for all $x\in \mathcal{M}$. The $r$-coverning number of $\mathcal{M}$ is the size of the smallest $r$-net. The {\em $r$-entropy number of $\mathcal{M}$}, $\mathsf{H}(\mathcal{M},r)$, is the natural log of the $r$-covering number. It is known that
\begin{equation}\label{eq:entropybound}\gamma^{(\alpha)}_2(\mathcal{M})\leq C\,\int_0^{\diam(\mathcal{M})}\,\sqrt{\mathsf{H}(\mathcal{M},r^{\frac{1}{\alpha}})}\,dr,\end{equation}
with $C>0$ is a universal constant. An important special case is when $\mathcal{M}\subset \re^d$ and $\dist$ is given by a norm, in which case the entropy integral bound is upper bounded by $\diam(\mathcal{M})^{\alpha}\sqrt{d}$ up to a universal constant. In particular, we obtain, 
\begin{equation}\label{eq:simpleboundgamma2}\gamma_{2}^{(\alpha)}(\mathcal{M})\leq C_{\alpha}\sqrt{d}\,\diam(\mathcal{M})^{\alpha}\end{equation}
with $C_{\alpha}>0$ only depends on $\alpha$. However, this bound can be very loose, as the next example shows. 

\begin{example}Let $\mathcal{M}\subset \re^d$ denote the standard simplex in $d$ dimensions, that is, the convex hull of the $d$ canonical basis vectors. Let $\dist$ denote the standard Euclidean metric. In this case, (\ref{eq:simpleboundgamma2}) bounds $\gamma_{2}^{(1)}(\mathcal{M}) = \mathcal{O}(\sqrt{d})$. By contrast, (\ref{eq:gamma2convexhull}) shows that $\gamma_{2}^{(1)}(\mathcal{M})$ is of the order of $\sqrt{\log d}$, which is sharp.\end{example}

\section{Setup and assumptions for main results}\label{sec:basic}

We now present the general setup and assumptions we will use in the analysis of SAA. 

\subsection{Ideal optimization versus SAA}\label{sub:basic}

\paragraph{Functions and sets.} As in the introduction, $\mathcal{I}$ is a finite set which will index the constraints of our problem. We use $0\not\in\mathcal{I}$ to index the objective function and set $\mathcal{I}_0:=\mathcal{I}\cup \{0\}$. 

We are given a set $Y\subset \re^d$ and functions $f_i:Y \to \re$, for $i\in \mathcal{I}_0$. We will also write $f:=f_0$. Given $\delta\in\re$, we define:
\[X_\delta:= \{x\in Y\,:\, \forall i\in\mathcal{I},\, f_i(x)\leq \delta\}.\]
We also write $X$ instead of $X_0$. Note that $X=X_{\delta}=Y$ for all $\delta>0$ when $\mathcal{I}=\emptyset$. The ``ideal"~optimization problem we consider is: 

\begin{eqnarray}
f^*:=\min_{x\in Y}&\quad & f(x) \label{problem:ideal}\\
\mbox{s.t.}&\quad &f_i(x)\leq 0,\, i\in\mathcal{I}.\nonumber
\end{eqnarray} 
In other words, the feasible set is $X$, the objective function is $f=f_0$ and the value of the problem is $f^*$. We will always assume implicitly that $X\neq \emptyset$. We let 
\[x^*\in{\rm arg\, min}_{x\in X}f(x)\mbox{ so that }f^*:=f(x^*).\]
In particular, we assume implicitly that our problem always has minimizers. We also use the symbols:
\[f^*_\delta:=\inf_{x\in X_\delta}f(x)\mbox{ and }{\rm gap}(\delta):=|f^*_\delta-f^*| \mbox{ (when $X_{\delta}\neq \emptyset$)}.\]
In case the above infimum is attained, we let 
\[x^*_{\delta}\in{\rm arg\, min}_{x\in X_\delta}f(x)\mbox{ so that }f^*_\delta:=f(x^*_\delta).\]

We will need some additional notation. We write $X_{\delta,{\rm act}(i)}$ for the subset of $X_\delta$ where constraint $i$ is active:
\[X_{\delta,{\rm act}(i)}:= \{x\in X_\delta\,:\, f_i(x) = \delta\}.\] We also define the set of points $x\in X_\delta$ that achieve $f(x)\leq f^*+\stheta$:
\[X_{\delta}^{*,\stheta}:= \{x\in X_\delta\,:\, f(x)\leq f^* + \stheta\}.\]
We set
\[X_{\delta}^{*,=\stheta}:= \{x\in X_\delta\,:\, f(x) =f^* + \stheta\}\]
and finally
\[X_{\delta,{\rm act}(i)}^{*,\stheta} := X^{*,\stheta}_\delta\cap X_{\delta,{\rm act(i)}}.\]
We emphasize that $\delta$ will be omitted from our notation when it is equal to zero. With few exceptions which are clear from context, we reserve the symbols $\delta,\eta$ for feasibility deviations and $\epsilon,\stheta$ for optimality deviations.

\paragraph{Randomness.} Let $(\Xi,\sigma(\Xi),\probn)$ denote a probability space. We write $\xi\sim \probn$ to denote a random element of $\Xi$ with law $\probn$.  In this paper,$\{\xi_k\}_{k=1}^N\subset \Xi$ is an i.i.d. random sample of size $N$ from the probability measure $\probn$. The $\xi_k$ are defined over a common probability space $(\Omega,\mathcal{A},\prob)$ that will be always kept implicit. $\widehat{\probn}$ denotes the empirical measure of the sample:
\[\widehat{\probn}:= \frac{1}{N}\sum_{k=1}^{N}\delta_{\xi_k}.\]
Given a measurable function $H:Y\times \Xi\to \re$ and $y$, we define:
\[\espn H(y,\cdot):= \int_{\Xi}H(y,\xi)\,\probn(d\xi)\mbox{ and }\widehat{\espn}\,H(y,\cdot) := \frac{1}{N}\sum_{k=1}^{N}H(y,\xi_k)\]
to denote the expectation and sample average (respectively) of $H(y,\cdot)$ with $y$ fixed. Our assumptions will be such that the integral over $\probn$ will always be well defined. 

\paragraph{Sample average approximation.} We are given measurable functions $F_i:Y \times \Xi\to \re$ for $i\in \mathcal{I}_0$. We assume that 
\begin{equation}\label{eq:SAAcondition} \forall i\in \mathcal{I}_0\,\forall x\in Y\,:\, \espn|F_i(x,\cdot)|<+\infty\mbox{ and }\espn F_i(x,\cdot) = f_i(x).\end{equation}
Write:\[\widehat{F}_i(x):=\widehat{\espn}F_i(x,\cdot) = \frac{1}{N}\sum_{k=1}^NF_i(x,\xi_k)\]
to denote the sample average of $F_i$. Formally, $\widehat{F}_i(x)$ is a function of $x$ and the sample, but we omit the sample from our notation. We sometimes write $\widehat{F}:=\widehat{F}_0$. The sample average approximation to problem (\ref{problem:ideal}) is:
\begin{eqnarray}
\widehat{F}^*:=\min_{x\in Y}&\quad & \widehat{F}(x) \label{problem:SAA}\\
\mbox{s.t.}&\quad &\widehat{F}_i(x)\leq 0,\, i\in\mathcal{I}.\nonumber
\end{eqnarray} 
 
Intuitively, the $\widehat{F}_i$ should give random approximations to the $f_i$ for large $N$, and optimization problems with the $\widehat{F}_i$ should be similar to the ``ideal" ~problems involving the $f_i$. Quantifying the extent to which this is true is the goal of this paper. We will need the following analogues of the notation introduced above:
\begin{eqnarray*}\widehat{X}&:=& \{x\in Y\,:\,\forall i\in\mathcal{I},\, \widehat{F}_i(x)\leq 0\};\\
\widehat{X}^{*,\stheta}&:=& \{x\in \widehat{X}\,:\, \widehat{F}(x)\leq \widehat{F}^* + \stheta\}.\\
\widehat{x}^* & \in & {\rm arg min}_{x\in \widehat X} \widehat{F}(x).\end{eqnarray*}

Again, we implicitly assume that the SAA always has solutions.

\subsection{Assumptions on the random functions} To state our general theorems, we will need some {\em probabilistic assumptions} on the random functions $F_i$. We start with a definition.

\begin{definition}[Good and great random variables] Given $(\sigma^2,\rho)\in\re_+$, a function $h:\Xi\to \re_+$ is said to be $(\sigma^2,\rho)$-good if $\espn\,h(\cdot)\leq \sigma^2$ and
\[\prob\left\{\widehat{\espn}h(\cdot)> 2\sigma^2\right\}\leq \rho.\]
Given $\sigma^2>0$, $p\geq 2$ and 
$\kappa_p>0$, we say that $h$ is $(\sigma^2,p,\kappa_p)$-great if  $\espn h(\cdot)\leq \sigma^2$ and in addition we have the $L^p$ norm bound:
\[\Lpnorm{h(\xi) - \espn h(\cdot)}\leq \kappa_p\sigma^2.\]\end{definition}
Any {\em fixed} integrable function $h\geq 0$ with $\espn h(\cdot)\leq \sigma^2$ is $(\sigma^2,\rho)$-good when $N$ is large enough due to the Law of Large Numbers. The point of our definition is to have finite-$N$ results. The next proposition says that great random variables satisfy a quantitative form of goodness.

\begin{proposition}[Proof in the Appendix] \label{prop:goodisgreath}If $h$ as above is 
$(\sigma^2,p,\kappa_p)$-great, it is also $(\sigma^2,\rho)$-good, with
\[\rho:=\left(\cbdg\kappa_p\sqrt{\frac{p}{N}}\right)^p\]
and $\cbdg$ is a universal constant.\end{proposition}

The kind of assumption we will make on the $F_i$ is described below. In what follows, $Z\subset Y$ is a subset of $Y$ containing $x^*$, $\sigma^2,\sigma_*^2,\rho>0$, $\alpha\in (0,1]$, $\kappa_p\geq 1$ and $p\geq 2$. Also, $\Vert\cdot\Vert$ is a norm over $\re^d$.

\begin{assumption}[$\parametersgoodness$-goodness over $Z$]\label{assump:goodFi} The functions $\{f_i\}_{i\in\mathcal{I}_0}$ and $ \{F_i\}_{i\in\mathcal{I}_0}$ are continuous in $x\in Y$. Moreover, 
\begin{enumerate}
\item The maps $\xi\mapsto (F_i(x^*,\xi) - f_i(x^*))^2\in \re_+$ are $(\sigma_*^2,\rho)$-good for each $i\in\mathcal{I}_0$;
\item For each map $F_i$ with $i\in\mathcal{I}_0$, there exists $\mathsf{L}_i:\Xi\to\re$ such that $\mathsf{L}_i^2$ is $(\sigma^2,\rho)$-good and:
\[\forall x,x'\in Z,\,\forall \xi \in\Xi\,:\, |F_i(x,\xi) - F_i(x',\xi)|\leq \mathsf{L}_i(\xi)\,\|x-x'\|^\alpha.\]\end{enumerate}
\end{assumption}

\begin{assumption}[$\parametersgreatness$-greatness over $Z$]\label{assump:greatFi} The functions $\{f_i\}_{i\in\mathcal{I}_0}$ and $ \{F_i\}_{i\in\mathcal{I}_0}$ are continuous in $x\in Y$. Moreover, 
\begin{enumerate}
\item The maps $\xi\mapsto (F_i(x^*,\xi) - f_i(x^*))^2\in \re_+$ are $(\sigma_*^2,p,\kappa_p)$-great for each $i\in\mathcal{I}_0$;
\item For each map $F_i$ with $i\in\mathcal{I}_0$, there exists $\mathsf{L}_i:\Xi\to\re$ such that $\mathsf{L}_i^2$ is $(\sigma^2,p,\kappa_p)$-great and:
\[\forall x,x'\in Z,\,\forall \xi \in\Xi\,:\, |F_i(x,\xi) - F_i(x',\xi)|\leq \mathsf{L}_i(\xi)\,\|x-x'\|^\alpha.\] \end{enumerate}
\end{assumption}

In our main results, we will make one of these two assumptions. For general problems, without a convexity assumption, we will take $Z=Y$. In convex settings, we will take potentially much smaller sets $X^{*,\stheta}_{\delta}\subset X$. Notice that each of the above assumptions implies:
\begin{equation}\label{eq:holdercont}\forall i\in\mathcal{I}_0,\,\forall x,x'\in Z\,:\, |f_i(x) - f_i(x')|\leq (\espn\mathsf{L}_i(\cdot))\,\|x-x'\|^{\alpha}\leq \sigma \,\|x-x'\|^{\alpha},\end{equation}
that is, the functions $f_i$ are $\alpha$-H\"{o}lder continuous over $Z$. 

We note the following simple consequence of Proposition \ref{prop:goodisgreath}.

\begin{proposition}[Great implies good; proof omitted] Assumption \ref{assump:greatFi} implies Assumption \ref{assump:goodFi} with the same set $Z$, the same parameters $\sigma^2,\sigma_*^2$, and 
\[\rho:=\left(\cbdg\kappa_p\sqrt{\frac{p}{N}}\right)^p.\]\end{proposition}

In particular, our assumptions may be satisfied with $\rho$ polynomially small in $N$, even if the random variables involved do not have light tails. 

\subsection{Assumptions on the geometry of the problem}\label{sub:geometry}

When there are constraints in expectation, the SAA will unavoidably have a different feasible set than the ideal problem. In this section, we present standard assumptions that allow us to bound the difference between the two sets. The first assumption is often used in the analysis of perturbations and algorithms for problems in Optimization and Variational Analysis \cite{pang,bauschke:borwein,iusem:jofre:thompson2015}. In what follows, $\|\cdot\|$ is a norm over $\re^d$ and $\dist$ is the corresponding set-to-point distance.

The first assumption is of Metric Regularity.
\begin{assumption}[Metric regular feasible (MRF) set]\label{assump:MRF} 
There exists $\mathfrak{c}>0$ such that for all $x\in Y$,
$$
\dist(x,X)\le \mathfrak{c}\sup_{i\in\mathcal{I}}f_i(x)_+.
$$

\end{assumption}
This assumption is trivially satisfied when $\mathcal{I}=\emptyset$.

MRF is related to standard constraint qualifications, e.g. the \emph{Slater constraint qualification} (SCQ) which ensures that $X$ has a strictly feasible point. For instance, Robinson \cite{robinson1975} proved that if the set $Y$ and the functions $f_i$ are convex, then, for some $\eta>0$,
\[X_{-2\eta}\neq \emptyset\Rightarrow  \mbox{Assumption \ref{assump:MRF}  holds with }\mathfrak{c}:=\frac{\diam(X)}{\eta}.
\]

The MRF condition is also true for a larger class of sets which are neither strictly feasible nor convex. One fundamental instance is of a polyhedron, as implied by Hoffmann's Lemma \cite{hoffman1952}. We remark here that, in Assumption \ref{assump:MRF}, we restrict our analysis for the case of ``Lipschitzian'' bounds. Our results can be easily extended to the case of ``H\"olderian'' bounds: for some $\beta>0$, $\dist(\cdot,X)\le\mathfrak{c}\sup_{i\in\mathcal{I}}\left[f_i(x)\right]_+^\beta$ (see Section 4.2 in \cite{pang}). In that case, MRF holds true for any compact nonconvex $X$ whose constraints are polynomial or real-analytic functions, a deep result implied by Lojasiewicz's inequality \cite{lojasiewicz1959}. We refer to Section 4.2 in \cite{pang} and references therein.

For convex problems, we will also consider a localized version of the Slater CQ condition. Here, we only require that the set $X^{*,\stheta}$ be bounded and has an ``interior point".

\begin{assumption}[Localized Slater CQ with convexity (LSCQ)]\label{assump:LSCQ} The set $Y$ is convex and closed, and the functions $\{f_i\}_{i\in\mathcal{I}_0}$ are continuous and convex. Moreover, there exist $\eta_*>0$ and $\stheta_*$ such that $X^{*,\stheta_*}$ is bounded and $X_{-\eta_*}^{*,\stheta_*}\neq \emptyset$ (that is, there exists $x\in Y$ with $f(x)\leq f^*+\stheta_*$ and $f_i(x)\leq -\eta_*$ for all $i\in\mathcal{I}$).
\end{assumption}

Boundedness of $X^{*,\stheta_*}$ may be guaranteed by usual assumptions. In that case, the Slater CQ, i.e., $X_{-\eta_*}\neq\emptyset$ for some 
$\eta_*>0$, implies Assumption \ref{assump:LSCQ} with $\stheta_* \ge \inf_{x\in X_{-\eta_*}}f(x)-f^* = {\rm gap}(-\eta_*)$. Assumption \ref{assump:LSCQ} allows us to control the complexity of $X^{*,\stheta}_\delta$ in terms of $X^{*,\stheta}$, for suitable $\stheta$ and $\delta$; see Lemma \ref{lem:robinson2} below for details.

We conclude this section by noting that in the next Sections \ref{sec:results}-\ref{sec:convexresults}, the functional $\gamma_2^{(\alpha)}$ is defined with respect to $\dist$, i.e., the set-to-point distance associated to the norm $\Vert\cdot\Vert$ over $\re^d$. See Assumptions \ref{assump:goodFi}, \ref{assump:greatFi} and \ref{assump:MRF}. 

\section{Main result for not-necessarily convex problems}\label{sec:results}

In this section we state formally and discuss our main result for SAA where we do not assume convexity. More precisely, we only make continuity and metric regularity assumptions on the functions we consider.  Theorem \ref{thm:generalrandomset} is closely related to previous results in the area. Our main contribution here is to obtain stronger bounds under light-tailedness assumptions, through the use of ``generic chaining"~and our novel concentration arguments. 

\begin{theorem}[General functions and sets; proof in \S \ref{proof:generalrandomset}]\label{thm:generalrandomset} Assume $Y$ is bounded. Additionally, make the assumption of $\parametersgoodness$-goodness over $Y$ (cf. Assumption \ref{assump:goodFi}). Given $t\geq 0$, define:
\[\widehat{r}_N(t):= \sigma\,\frac{\cgammaconc\,\gamma^{(\alpha)}_2(Y) + \cdiamconc\,\diam^{\alpha}(Y)\sqrt{1+ \log (2|\mathcal{I}|+2) + t}}{\sqrt{N}}.\]
Also define:
\begin{eqnarray*}\widehat{\delta}_N(t) &:=& \left\{\begin{array}{ll} 
\sigma_*\sqrt{\frac{6(1+\log (2|\mathcal{I}|+2) + t }{N}} +\widehat{r}_N(t), & \mathcal{I}\neq \emptyset,\\ 0, & \mathcal{I}=\emptyset. \end{array}\right.
\end{eqnarray*}
Let ${\rm Good}_{\rm Thm. \ref{thm:generalrandomset}}(\epsilon_0,t)$ denote the event where the following properties hold:
\begin{enumerate}
\item[{\bf (a)}] $\widehat{X}\subset X_{\widehat{\delta}_N(t)}$, that is, feasible points of the SAA violate ideal constraints by at most $\widehat{\delta}_N(t)$;
\item[{\bf (b)}] \[\widehat{X}^{*,\epsilon_0}\subset  X^{*,\epsilon_0 + 2\widehat{r}_N(t) + {\rm gap}(-\widehat{\delta}_N(t))}_{\widehat{\delta}_N(t)},\] that is, for any $x\in \widehat{X}$ with $\widehat{F}(x)\leq \widehat{F}^*+\epsilon_0$, we have $f(x)\leq f^*+\epsilon_0 + {\rm gap}(-\widehat{\delta}_N(t)) + 2\widehat{r}_N(t)$ and $\max_{i\in\mathcal{I}}f_i(x)\leq \widehat{\delta}_N(t)$ (recall that ${\rm gap}(\delta):=|f^*_\delta-f^*|$, cf. \S \ref{sub:basic});
\item[{\bf (c)}] \[|\widehat{F}^* - f^*| \leq \widehat{\delta}_N(t) +
\max\left\{2\,\widehat{r}_N(t) + {\rm gap}(-\widehat{\delta}_N(t)),\gap(\widehat{\delta}_N(t))\right\}.\] 
\end{enumerate}
Then $\prob({\rm Good}_{\rm Thm.\ref{thm:generalrandomset}}(\epsilon_0,t))\geq 1 - e^{-t}-2(|\mathcal{I}|+1)\rho$.
If we assume $\parametersgreatness$-greatness over $Y$ (cf. Assumption \ref{assump:greatFi}) instead of $\parametersgoodness$-goodness, then one may take $\rho = (\cbdg\kappa_p\sqrt{p/N})^p$ above. Finally, if we additionally make the metric regularity assumption (Assumption \ref{assump:MRF}), we have the following inequality whenever ${\rm Good}_{\rm Thm. \ref{thm:generalrandomset}}(\epsilon_0,t)$ occurs:
 \[\forall x\in \widehat{X}, \dist(X,x)\leq  \mathfrak{c}\,\widehat{\delta}_N(t).\] 
 \end{theorem}

Let us parse this theorem. The error parameter $\widehat{r}_N(t)+{\rm gap}(-\widehat{\delta}_N(t))$ controls how good SAA solutions are for the original problem. The related parameter 
$\widehat{\delta}_N(t)+\widehat{r}_N(t)+\max\{{\rm gap}(-\widehat{\delta}_N(t)),\gap(\widehat{\delta}_N(t))\}$ bounds the difference between values of the SAA and the ideal problem. Finally, $\widehat{\delta}_N(t)$ controls just how much SAA feasible points violate the constraints of the original problem, and also (under Assumption \ref{assump:MRF}) how far feasible points of the SAA are from the ideal feasible set. Note that if Assumption \ref{assump:MRF} holds on  the set $X_{-\delta_*}$ for some small $\delta_*>0$, $\max\{\gap(\delta),\gap(-\delta)\}$ is of the order of $\mathfrak{c}\delta^\alpha$ for any $\delta\in[0,\delta_*]$.

The main features of these parameters is their dependence on the sample size $N$,  the geometry of the problem and the desired probability level. The dependence on $N$ is always of the form $N^{-1/2}$, in contrast with previous analyses of SAA not requiring light tails \cite{kancova:omelchenko2015}. The {\em geometry} of the set $Y$ comes into play via the diameter of $Y$ and the Gaussian complexity parameter $\gamma_2^{(\alpha)}(Y)$. These parameters are optimal for controlling fluctuations of Gaussian processes, and we show that they may still be used in heavier-tailed settings. Finally, the error bounds depend in a sub-Gaussian fashion on the desired probability level $e^{-t}$, at least when $(2|\mathcal{I}|+1)\rho\ll e^{-t}.$ In that connection, we note $\rho$ decays polynomially with $N$ under the $\parametersgreatness$-greatness assumption, if $p$ and $\kappa_p$ is treated as a constant. Therefore, our Theorem \ref{thm:generalrandomset} does give sub-Gaussian-type error probabilities if the number of constraints satisfies $|\mathcal{I}|\leq N^{p/2-c}$ (with $c>0$) and $t \leq c \log N$. We expect this to be usually the case in applications. Still, we observe that our assumptions for Theorem \ref{thm:generalrandomset} are somewhat limiting, as they do {\em not} allow for unbounded feasible sets (for example).

\begin{remark}\label{rem:thm:gen:slack}
A natural question is if there are advantages in considering the SAA feasible set 
$\widehat{X}=\{x\in Y:\widehat F_i(x)\le\delta,\forall i\in\mathcal{I}\}$ with a positive slack 
$\delta>0$. A corollary of the proof of Theorem \ref{thm:generalrandomset} is that by choosing 
$\delta:=\mathcal{O}(\widehat{\delta}_N(t))$, we can remove $\gap(-\widehat{\delta}_N(t))$ in the bounds of item \textbf{(b)} and \textbf{(c)} above under essentially the same assumptions. Of course, this is of theoretical interest only as the constants in the rate 
$\widehat{\delta}_N(t)$ are typically unknown. 
\end{remark}

\section{Main result in the convex case}\label{sec:convexresults}

We now consider a situation where Theorem \ref{thm:generalrandomset} can be improved upon. By assuming that the set $Y$ and the functions $F_i$ are convex and the feasible set satisfies a localized Slater-type condition (Assumption \ref{assump:LSCQ}), we will see that we can obtain a stronger result, Theorem \ref{thm:convexrandomset} below. Before we present it, we first discuss some geometrical aspects of the problem, which will explain the somewhat convoluted form of the theorem.

\subsection{A preliminary discussion}\label{ss:prelim:discussion}

Throughout this section, we make Assumption \ref{assump:LSCQ} that the set $X^{*,\vartheta_*}_{-\eta_*}\neq \emptyset$. This implies that there exists a point $x\in X$ in the feasible set of the ideal problem that satisfies the following properties:

\begin{enumerate}
    \item $\max_{i\in\mathcal{I}}f_i(x)\leq -\eta_*$ for some $\eta_*>0$, that is, all constraints are ``far"~from being active on $x$;
    \item $f(x)\leq f^*+\vartheta_*$ for some $\vartheta_*\geq 0$, that is, $x$ is a near optimizer of the ideal problem.
\end{enumerate}

As noted in the discussion after Assumption \ref{assump:LSCQ}, we can assume that $\vartheta_* \geq {\rm gap}(-\eta_*)$. We are especially interested in situations where $X^{*,\vartheta_*}_{\eta_*}$ is bounded; this is the case for instance if $f$ has a unique minimizer $x_*$, satisfies a growth condition $f(x)-f^*\geq c\,\min\{|x-x_*|^\beta,1\}$, for some constants $\beta,c>0$ and $\eta_*$ small enough. Notice that $X^{*,\vartheta_*}_{\eta_*}$ can be bounded while the whole set $X$ is unbounded. 

We now consider the role of convexity. Recall that $\widehat{x}$ is a solution to the SAA. Given $\delta\in (0,\eta_*]$ and $\epsilon\in [{\rm gap}(-\delta),\vartheta_*]$. Say that $\widehat{x}$ is $(\delta,\epsilon)$-{\em good} if:
\begin{itemize}
    \item no constraint of the original problem is violated by more than $\delta$: \[\max_{i\in\mathcal{I}}f_i(\widehat{x})\leq \delta;\]
    \item the objective function at $\widehat{x}$ satisfies $f(\widehat{x})\leq f^*+\epsilon$.
\end{itemize}
We say $\widehat{x}$ is $(\delta,\epsilon)$-{\em bad} if it is not $(\delta,\epsilon)$-{good}. What could cause $\widehat{x}$ to be bad? Proposition \ref{lem:convexlocalized}, a deterministic result, shows that, if $\widehat{x}$ is $(\delta,\epsilon)$-bad, then there exists a point $x\in X^{*,\epsilon}_\delta$ where $\widehat{F}_i(x)-f_i(x)$ is ``large" for some $i\in\mathcal{I}\cup \{0\}$ (recall that $i=0$ corresponds to the objective function). That is, if the SAA solution is bad, this is due to a failure of concentration of the SAA functions around their ideal counterparts. Most importantly, this failure must happen in the set $X^{*,\epsilon}_\delta$, which will often be much smaller than $X$ (it is at most as large as $X^{*,\vartheta_*}_{\eta_*}$). This is what we mean by {\em localization}: failure of the SAA manifests itself at ``small scales". 

As a second step, we further analyze the set $X^{*,\epsilon}_\delta$. It will follow from Lemma \ref{lem:robinson2} that 
\[X^{*,\epsilon}_\delta\subset 2X^{*,\epsilon}-x_{-\delta},\]
for some point $x_{-\delta}\in X^{*,\epsilon}_{-\delta}$. This means that $X^{*,\epsilon}_\delta$ is contained in a homothetic copy of $X^{*,\epsilon}$. As noted above, if $f$ satisfies a growth assumption, the diameter of $X^{*,\epsilon}$ goes to $0$ as $\epsilon\searrow 0$. In particular, this will mean that $X^{*,\epsilon}_\delta$ is also small.

The upshot of our discussion so far is this. Suppose we can suitably guarantee that, with high probability, we have that the sample averages $\widehat{F}_{i}(x)$ are uniformly close to $f_i(x)$ for all $i\in\mathcal{I}\cap\{0\}$, for all $x\in X^{*,\epsilon}_\delta\subset 2X^{*,\epsilon}-x_{-\delta}$. Then it follows that the $\widehat{x}$ is $(\delta,\epsilon)$-good. 

How does one choose $\delta$ and $\epsilon$ that are as small as possible, while ensuring that $\widehat{x}$ is $(\delta,\epsilon)$-good with high probability? As it turns out, this is somewhat tricky. To a first approximation, we should expect that:  
\begin{align}
\sup_{x\in X^{*,\epsilon}_{\delta}}|\widehat{F}_{i}(x) - f_i(x) - \widehat{F}_i(x_*) + f(x_*)| \approx \sigma(\epsilon,\delta)\frac{\gamma^{(\alpha)}_2(X^{*,\epsilon})}{\sqrt{N}},\label{equation:local:complexity:heuristic}
\end{align}
where $\sigma(\epsilon,\delta)$ is a term pertaining to the Lipschitz or H\"{o}lder constants of the functions $\widehat{F}_i$ over the set $X^{*,\epsilon}_\delta$. The conditions we need are that these and other random quantities are smaller than both $\delta$ and $\epsilon$, so that the ``noise"~terms do not overwhelm the ``signal" in the SAA. Such difficulties also appear in the literature on localization in Statistics and Machine Learning \cite{koltchinskii2006,mendelson2017}, and lead to somewhat convoluted statements. This literature however assume fixed constraints. Our setting study localization with random constraints and one has to account for the fact that $\delta$ and $\epsilon$ are \emph{coupled} via $\epsilon\in [{\rm gap}(-\delta),\vartheta_*]$. In any case, the parameter choices in Theorem \ref{thm:convexrandomset} will be derived from variants of the above reasoning. In most typical situations, one has available upper bounds on the ``local complexities'' defined in the right hand side of \eqref{equation:local:complexity:heuristic}. See e.g. \eqref{eq:simpleboundgamma2}. In this case, the above reasoning leads to solving a ``fixed-point'' equation in $(\delta,\epsilon)$. While difficult to solve in general, sufficient upper bounds can be obtained by solving inequalities in $(\delta,\epsilon)$. We exemplify this reasoning in Section \ref{ss:particular:case} and \ref{sec:simpleresults}.

\subsection{The theorem}
 
 We can now state the main result of this section. 

\begin{theorem}[Convex sets and functions; proof in \S \ref{proof:convexrandomset}]\label{thm:convexrandomset} Make Assumption \ref{assump:LSCQ} with constants $\eta_*,\stheta_*$. Also assume $\parametersgoodnessplus$-goodness over the set $Z=X^{*,\stheta}_{\delta}$ for every choice of $(\stheta,\delta)\in [0,\stheta_*]\times [0,\eta_*]$ (cf. Assumption \ref{assump:goodFi}), where $\sigma^2(\stheta,\delta)$ depends continuously on $\stheta$ and $\delta$ (note that $\sigma^2(\stheta,\delta)$ depends on $(\stheta,\delta)$ but the other parameters in Assumption \ref{assump:goodFi} are fixed). 

Fix parameter $t\geq 0$. For every $0<\epsilon\leq \stheta_*$ and $0<\delta<\eta_*$ satisfying $\epsilon+{\rm gap}(-\delta)\leq \stheta_*$, set:
\begin{eqnarray*}\widehat{w}_N(t;\delta;\epsilon)&:=&  \,\sigma(\epsilon+{\rm gap}(-\delta);\delta)\left\{4\sqrt{3}\frac{\gamma^{(\alpha)}_2(X^{*,\epsilon+{\rm gap}(-\delta)})}{\sqrt{N}}\right. \\ & & \left. + 6\sqrt{3}\frac{\diam^{\alpha}(X^{*,\epsilon+{\rm gap}(-\delta)})\sqrt{1 + \log(2|\mathcal{I}|+2)+t}}{\sqrt{N}}\right\},\end{eqnarray*} 
For $(\epsilon,\delta)$ as above, we define parameters $\check{\delta}(t;\epsilon)$ and $\check{w}(t;\epsilon)$ as follows. 
\begin{enumerate}\item If $\mathcal{I}=\emptyset$ (there are no constraints in expectation), then $\check{\delta}(t;\epsilon):=0$ and $\check{w}(t;\epsilon) = \widehat{w}_N(t;0;\epsilon)$.
\item Otherwise, assume that \[S_{N,\eta_*}(t;\epsilon):=\left\{\delta\in(0,\eta_*) \,:\, \begin{array}{l}\epsilon + {\rm gap}(-\delta) \leq \stheta_*,\\ \widehat{w}_N(t;\delta;\epsilon)+ \sigma_*\sqrt{\frac{6(1+\log(2|\mathcal{I}|+2)  + t)}{N}}< \delta\end{array}\right\}\]
is nonempty, and define \[\check{\delta}(t;\epsilon):=\inf S_{N,\eta_*}(t;\epsilon)\mbox{ and }\check{w}(t;\epsilon) := \widehat{w}_N(t;\check{\delta}(t;\epsilon);\epsilon).\]\end{enumerate}

Now, fix $\epsilon_0\in[0,\vartheta_*)$ and assume the set 
\[R_{N,\eta_*}(t;\epsilon_0):=\{\epsilon \in (\epsilon_0,\stheta_*]\,:\, \epsilon > \epsilon_0 + {\rm gap}(-\check{\delta}(t;\epsilon)) +  2\check{w}(t;\epsilon)\},\]
is nonempty so that 
$$
\check{r}(t;\epsilon_0):=\inf R_{N,\eta_*}(t;\epsilon_0)
$$ 
is well defined. Also set 
$$\check{\delta}(t):=\lim_{\epsilon\searrow \check{r}(t;\epsilon_0)}\check{\delta}(t;\epsilon).
$$

Now define ${\rm Good}_{\rm Thm. \ref{thm:convexrandomset}}(t,\epsilon_0)$ as the event where the following properties all hold.
\begin{enumerate}
\item[{\bf (a)}] \[\widehat{X}^{*,\epsilon_0}\subset X_{\check{\delta}(t)}^{*,\check{r}(t;\epsilon_0)};\]
that is, all $x\in \widehat{X}$ with $\widehat{F}(x)\leq \widehat{F}+\epsilon_0$ also satisfy $f(x)\leq f^*+\check{r}(t;\epsilon_0)$ and $\max_{i\in\mathcal{I}}f_i(x)\leq \check{\delta}(t)$;
\item[{\bf (b)}] the values of the SAA and the ideal problem satisfy:
\[|\widehat{F}^* - f^*|\leq \sigma_*\sqrt{\frac{6(1+\log(2|\mathcal{I}|+2)  + t)}{N}} + \frac{\check{r}(t;\epsilon_0)}{2}
+\max\{\check{r}(t;\epsilon_0),\gap(\check\delta(t))\}.
\]

\item[{\bf (c)}] for all $x\in \widehat{X}^{*,\epsilon_0}$, \[\dist(x,X)\leq \min\left\{\frac{\diam(X^{*,\stheta_*})\,\check{\delta}(t)}{\eta_*},2\diam(X^{*,\check{r}(t;\epsilon_0)})\right\}.\]
\end{enumerate}
Then \[\prob({\rm Good}_{\rm Thm. \ref{thm:convexrandomset}}(t,\epsilon_0))\geq 1-e^{-t}-2(|\mathcal{I}|+1)\rho.\] If we assume instead $\parametersgreatnessplus$-greatness of the functions $F_i$  (cf. Assumption \ref{assump:greatFi}) instead of $\parametersgoodnessplus$-goodness, then one may take $\rho = (\cbdg\kappa_p\sqrt{p/N})^p$ above.
\end{theorem}

The comments we made on Theorem \ref{thm:generalrandomset} on probabilities of error still apply. However, the statement of  Theorem \ref{thm:convexrandomset} is more convoluted. To begin with, the sets  $S_{N,\eta_*}(\cdot;\cdot)$ and $R_{N,\eta_*}(\cdot;\cdot)$ essentially constrain the choices of $\delta$ and $\epsilon$ so that (in the parlance of the preliminary discussion) the ``signal"~terms are always larger than the stochastic ``noise"~in the SAA. Nonemptyness of these sets, which is assumed in Theorem \ref{thm:convexrandomset}, is a consequence of a lower bound on the sample size $N$. The infima taken over these sets in Theorem \ref{thm:convexrandomset} correspond to trying to find the smallest possible $\delta$ and $\epsilon$ to which our reasoning applies, which are given by $\check{\delta}(t)$ and $\check{r}(t;\epsilon_0)$ (respectively). As with the obtained rates, a localized lower bound on $N$ can be obtained by using the control on the quantities  $\gap(-\delta)$, $\sigma(\epsilon+\gap(-\delta);\delta)$, $\gamma_2^{(\alpha)}(X^{*,\epsilon+\gap(-\delta)})$ and $\diam(X^{*,\epsilon+\gap(-\delta)})$ and solving the inequalities defining $S_{N,\eta_*}(\cdot;\cdot)$ and $R_{N,\eta_*}(\cdot;\cdot)$. This will be exemplified in Sections \ref{ss:particular:case} and \ref{sec:simpleresults}. 

Let us now discuss $\check{\delta}(t)$. Basically, this parameter controls fluctuations in the random constraints of the SAA. On the one hand, if we assume  
\begin{equation}\label{eq:interiorconvexrandomset}\begin{array}{c}\mathcal{I}=\emptyset\\ \mbox{{\bf or}}\\ \max_{i\in\mathcal{I}}f_i(x^*)\leq -\eta_0,\mbox{ where }\eta_0:=\sigma^*\sqrt{\frac{6(1+\log(4|\mathcal{I}|+2) + t)}{N}},\end{array}\end{equation}
we may then take take $\eta_*:=\eta_0$ and note $x^*\in X^{*,\stheta_*}_{-\eta_*}$, so that ${\rm gap}(-\delta)=0$ for all $0\leq \delta\leq \eta_*$. Intuitively, what this means is that $x^*$ satisfies the constraints with enough slack that it is nearly certain to be feasible for the SAA, in which case the random constraints do not matter much. 

Now assume (\ref{eq:interiorconvexrandomset}) does {\em not} hold. This means that there are random constraints and $x^*$ is on or near the boundary of the feasible set of the ideal problem. In particular, it may not be feasible for the SAA. However, the {\em existence} of a point $x_{-\eta_*}\in X^{*,\stheta_*}_{-\eta_*}$ gives stability results. Lemma \ref{lem:robinson2} below implies:
\begin{equation}
\label{eq:boundaryptdelta}
\forall 0\leq \delta<\eta_*\,:\,\,\exists x_{-\delta}\in X^{*,\stheta_*}_{-\delta}\,:\, \|x_{-\delta} - x^*\|\leq \frac{\diam(X^{*,\stheta_*})\,\delta}{\eta_*-\delta}.
\end{equation}
The goodness assumption over $Z=X^{*,\stheta_*}_{-\delta}$ in Theorem \ref{thm:convexrandomset}, i.e.,  (\ref{eq:holdercont}) with $\sigma:=\sigma(\vartheta_*,-\delta)$, gives 
\begin{equation}
\label{eq:boundaryptdelta:gap}
{\rm gap}(-\delta)\leq f(x_{-\delta}) - f(x^*)\leq \sigma\,{ \diam}^\alpha(X^{*,\stheta_*})\,\left(\frac{\delta}{\eta_*-\delta}\right)^{\alpha}.
\end{equation}
One can then use this bound on ${\rm gap}(-\delta)$ and the regularity conditions of the objective function and constraints to obtain upper bounds on $\check{r}(t;\epsilon_0)$ and $\check{\delta}(t)$. In general, this may lead to  bounds that can be significantly larger than when (\ref{eq:interiorconvexrandomset}) holds. We will see in \S \ref{sub:applyconvexrandomset} (especially in Remark \ref{rem:onedimsimple}) that such larger bounds are unavoidable in general even for simple metric projection problems.  

\begin{remark}\label{rem:thm:convex:slack}
As in Remark \ref{rem:thm:gen:slack}, there are advantages in considering 
$\widehat{X}=\{x\in Y:\widehat F_i(x)\le-\delta,\forall i\in\mathcal{I}\}$ with a slack 
$\delta>0$. A corollary of the proof of Theorem \ref{thm:convexrandomset} is that by taking
$\delta:=\mathcal{O}(\check{\delta}(t))$ with similar assumptions, it is possible to improve item \textbf{(a)} to 
$\widehat X^{*,\epsilon_0}\subset X^{*,\check r(t;\epsilon_0)}$ and remove $\gap(\check\delta(t))$ in the bound of item \textbf{(b)}. Again,  tuning $\delta$ to the order of $\check{\delta}(t)$ is of theoretical interest only as the latter is typically unknown. 
\end{remark}

\subsection{Two instructive particular cases}\label{ss:particular:case}

We finish this section with an application of the general Theorem \ref{thm:convexrandomset} when the solution set satisfies a local regularity condition. The purpose here is to further clarify the usefulness of ``localization'' (as discussed in Section \ref{ss:prelim:discussion}) in a typical setting in stochastic convex optimization. 

\begin{assumption}[Locally regular solution set]\label{assump:LRSS} Suppose that there exist $\mathfrak{c}>0$, $\kappa\in(0,1]$ and $\vartheta_*>0$ such that for all $\epsilon\in(0,\vartheta_*]$,
\begin{align*}
    \diam(X^{*,\epsilon})\le \mathfrak{c}\epsilon^\kappa+\diam(X^*).
\end{align*}
\end{assumption}

In the following,
\begin{align*}
\mathfrak{C}_\alpha:=
\sup_{\epsilon\in[0,\vartheta_*]}\left(\frac{\gamma_2^{(\alpha)}(X^{*,\epsilon})}{\diam^{\alpha}(X^{*,\epsilon})}\right)^2.
\end{align*}
A conservative bound is  $\mathfrak{C}_\alpha\le C_\alpha d$ for a constant $C_\alpha$ that depends only on $\alpha$. See Section \ref{ss:complexity:parameter}. 

Assumption \ref{assump:LRSS} deserves some discussion. One typical instance of Assumption \ref{assump:LRSS} is when $f$ is strongly convex on $Y$ (in this case, $\kappa=1/2$). More generally, Assumption \ref{assump:LRSS} is implied when $f$ is locally strongly convex on an open neighbourhood of $X^*$.\footnote{When $\mathcal{I}\neq\emptyset$, we assume without too much loss in generality in Assumption \ref{assump:LRSS} that $\epsilon\in(0,\vartheta^*]$ with $\vartheta^*$ as in Assumption \ref{assump:LSCQ}. For instance, in case $f$ is \emph{locally} strongly convex on a neighbourhood $U$ of $X^*$ and Assumption \ref{assump:LSCQ} holds, the existence of a $x\in X_{-\eta_*}^{*,\vartheta_*}\cap U$ is a mild requirement.} Other important instance when Assumption \ref{assump:LRSS} holds is when the problem has (local) \emph{weak sharp minima} \cite{1993burke:ferris,2005burke:deng} (in this case with $\kappa=1$).
 
We present two results, one when the feasible set $X$ is fixed and the second when $X$ has random constraints. 

\begin{proposition}[Fixed feasible set]\label{prop:LRSS}
Assume that $\mathcal{I}=\emptyset$ (that is, there are no random constraints). Grant Assumption \ref{assump:LRSS} with constants $\mathfrak{c},\kappa$ and $\vartheta_*$. The set $Y$ is convex and closed, and the functions $\{f_i\}_{i\in\mathcal{I}_0}$ are continuous and convex. Assume $(\sigma_*^2,\sigma^2,\alpha,\rho)$-goodness over the set $Z=X^{*,\stheta_*}$ (cf. Assumption \ref{assump:goodFi}). Let $t>0$ and $\epsilon_0\in[0,\vartheta_*/2)$.

Define 
$
\phi_N(t):=\sqrt{(\mathfrak{C}_\alpha+t)/N},
$
and, suppose $N$ large enough so that, for an absolute constant $C>0$,
\begin{align}
\epsilon_0+2C\sigma\left(\mathfrak{c}^\alpha\vartheta_*^{\kappa\alpha}+\diam^\alpha(X^*)\right)\cdot\phi_N(t)\le\vartheta_*/2.\label{prop:LRSS:N}
\end{align}
Define
\begin{align*}
\mathfrak{R}_N:=
\begin{cases}
2\epsilon_0+4C\sigma\diam^\alpha(X^*)\cdot\phi_N(t),&\mbox{ if $\alpha\kappa=1$},\\
\max\left\{2\epsilon_0+4C\sigma\diam^\alpha(X^*)\cdot\phi_N(t),\left[4C\mathfrak{c}^{\alpha}\sigma\phi_N(t)\right]^{\frac{1}{1-\alpha\kappa}}\right\},&\mbox{ if $\alpha\kappa\in(0,1)$}.
\end{cases}
\end{align*}

Finally, define ${\rm Good}_{\rm Prop. \ref{prop:LRSS}}(t,\epsilon_0)$ as the event where the following properties all hold.
\begin{enumerate}
\item[{\bf (a)}] \[\widehat{X}^{*,\epsilon_0}\subset X^{*,\mathfrak{R}_N};\]
that is, all $x\in \widehat{X}$ with $\widehat{F}(x)\leq \widehat{F}+\epsilon_0$ also satisfy $f(x)\leq f^*+\mathfrak{R}_N$ and $x\in X$;
\item[{\bf (b)}] the values of the SAA and the ideal problem satisfy:
\[|\widehat{F}^* - f^*|\leq \sigma_*\sqrt{\frac{6(1+\log 2  + t)}{N}} + \frac{3}{2}\mathfrak{R}_N.
\]
\end{enumerate}
Then \[\prob({\rm Good}_{\rm Prop. \ref{prop:LRSS}}(t,\epsilon_0))\geq 1-e^{-t}-2\rho.\] If we assume instead $(\sigma_*^2,\sigma^2,\alpha,p,\kappa_p)$-greatness of the functions $F_i$   (cf. Assumption \ref{assump:greatFi}) instead of $(\sigma_*^2,\sigma^2,\alpha,\rho)$-goodness, then one may take $\rho = (\cbdg\kappa_p\sqrt{p/N})^p$ above.
\end{proposition}
\begin{proof}
We only present a proof sketch. Denote $\phi:=\phi_N(t)$. As there is no random constraints, $\gap(-\delta)\equiv0$ and we may set $\check\delta(t,\epsilon)\equiv0$. In particular, $\check\delta(t)=0$ and, by the definitions of $\mathfrak{C}_\alpha$, $\phi_N(t)$ and $\check w(t;\epsilon):=\widehat w_N(t;0;\epsilon)$,
\begin{align*}
\check w(t;\epsilon)=C\sigma\phi\diam^\alpha(X^{*,\epsilon})
\le C\sigma\phi\mathfrak{c}^\alpha\epsilon^{\alpha\kappa}
+C\sigma\phi\diam^\alpha(X^*),
\end{align*}
for some constant $C>0$ and for all $\epsilon\in(0,\vartheta_*]$ by Assumption \ref{assump:LRSS}. Let us define $\mathsf{A}_N(\epsilon):=C\sigma\mathfrak{c}^\alpha\phi\epsilon^{\kappa\alpha}$ and $\mathsf{B}_N:=C\sigma\diam^\alpha(X^*)\phi$. 

Recall that $\check r(t;\epsilon_0)$ is the infimum over 
$\epsilon$ with constraints $\epsilon_0+2\check w(t;\epsilon)<\epsilon\le\vartheta_*.
$
An upper bound on $\check r(t;\epsilon_0)$ is obtained by considering the infimum over the smaller set $R$ defined by $\epsilon$ such that
$
\epsilon_0+2\mathsf{B}_N+2\mathsf{A}_N(\epsilon)
<\epsilon\le\vartheta_*
$. 
If \eqref{prop:LRSS:N} holds then $R\neq\emptyset$. Suppose first $\kappa\alpha\in(0,1)$. A simple calculation yields $\check r(t;\epsilon_0)\le\max\{2\epsilon_0+4\mathsf{B}_N,(4C\mathfrak{c}^{\alpha}\sigma\phi)^{\frac{1}{1-\alpha\kappa}}\}$. Suppose now $\alpha\kappa=1$. Again using \eqref{prop:LRSS:N}, a simple calculation shows that $\check r(t;\epsilon_0)\le2\epsilon_0+4\mathsf{B}_N$. This finishes the proof.
\qed
\end{proof}

For instance, in case of quadratic growth ($\kappa=1/2$) and Lipschitz continuity ($\alpha=1$), one has the optimality slackness $\mathfrak{R}_N$ of the order $\epsilon_0+\mathfrak{c}\sigma^2(\mathfrak{C}_\alpha+t)/N$. A notable feature of Proposition \ref{prop:LRSS} is that the ``rate'' $\mathfrak{R}_N$ on the sample size $N$ is \emph{independent} of $\diam(X)$. In particular, it allows unbounded $X$. This is in large contrast with the bounds obtained in Theorem \ref{thm:generalrandomset} in the general non-convex case. ``Localization'', implied by convexity, is the technique allowing for such sharper rates. When $\alpha\kappa=1$ and the solution is unique, so that $\diam(X^*)=0$, Proposition \ref{prop:LRSS} implies that for large enough $N$, the SAA solution is an \emph{exact} solution of the original problem with high probability. For convex piece-wise linear programs, this was been observed in \cite{shapiro:tito2000}.

We now consider the case of random constraints. 
\begin{proposition}[Random feasible set]\label{prop:LRSS:random:set}
Make Assumption \ref{assump:LSCQ} with constants $\eta_*,\vartheta_*$. Also assume $(\sigma_*^2,\sigma^2,\alpha,\rho)$-goodness over the set $Z=X^{*,\stheta_*}_{\delta_*}$  (cf. Assumption \ref{assump:goodFi}). Grant Assumption \ref{assump:LRSS} with constants $\mathfrak{c},\kappa$ and $\vartheta_*$. Assume (for simplicity) that $0<\alpha\kappa<1$. Let $t>0$ and $0\le\epsilon_0<\min\{\vartheta_*/4,1/2\}$. 

Then there is constant $\mathfrak{C}>1$ depending only on $\eta_*$, $\vartheta_*$, $\sigma$, $\sigma_*$, $\alpha$, $\kappa$, $\mathfrak{c}$ and $\diam(X^*)$ for which the following statement holds.
Let 
$$
\phi_N(t):=\sqrt{(\mathfrak{C}_\alpha+\log|\mathcal{I}|+t)/N},
$$
and assume that $N$ is large enough so that $\mathfrak{C}\phi_N(t)\le1$. Let \begin{align*}
\mathfrak{R}_N:=\max\left\{
2\epsilon_0,\mathfrak{C}\phi_N^\alpha(t)
\right\},\quad\mbox{and}\quad
\mathfrak{D}_N:=
\mathfrak{C}\phi_N(t)+\mathfrak{C}\phi_N(t)\max\left\{
\epsilon_0^{\alpha\kappa},\phi_N^{\alpha^2\kappa}(t)
\right\}.
\end{align*}
Finally, define ${\rm Good}_{\rm Prop. \ref{prop:LRSS:random:set}}(t,\epsilon_0)$ as the event where the following properties all hold.
\begin{enumerate}
\item[{\bf (a)}] \[\widehat{X}^{*,\epsilon_0}\subset X_{\mathfrak{D}_N}^{*,\mathfrak{R}_N};\]
that is, all $x\in \widehat{X}$ with $\widehat{F}(x)\leq \widehat{F}+\epsilon_0$ also satisfy $f(x)\leq f^*+\mathfrak{R}_N$ and $\max_{i\in\mathcal{I}}f_i(x)\leq \mathfrak{D}_N$;
\item[{\bf (b)}] the values of the SAA and the ideal problem satisfy:
\[|\widehat{F}^* - f^*|\leq \sigma_*\sqrt{\frac{6(1+\log(2|\mathcal{I}|+2)  + t)}{N}} + \frac{\mathfrak{R}_N}{2}
+\max\{\mathfrak{R}_N,\gap(\mathfrak{D}_N)\}.
\]

\item[{\bf (c)}] for all $x\in \widehat{X}^{*,\epsilon_0}$, \[\dist(x,X)\leq \min\left\{\frac{\diam(X^{*,\stheta_*})\,\mathfrak{D}_N}{\eta_*},2\diam(X^{*,\mathfrak{R}_N})\right\}.\]
\end{enumerate}
Then \[\prob({\rm Good}_{\rm Prop. \ref{prop:LRSS:random:set}}(t,\epsilon_0))\geq 1-e^{-t}-2(|\mathcal{I}|+1)\rho.\] If we assume instead $(\sigma_*^2,\sigma^2,\alpha,p,\kappa_p)$-greatness of the functions $F_i$ over $Y$  (cf. Assumption \ref{assump:greatFi}) instead of $(\sigma_*^2,\sigma^2,\alpha,\rho)$-goodness, then one may take $\rho = (\cbdg\kappa_p\sqrt{p/N})^p$ above.
\end{proposition}
\begin{proof}
Let $0<\epsilon\leq \min\{\stheta_*/2,1\}$ and $0<\delta\le\min\{\eta_*/2,1\}$ with $\epsilon_1:=\epsilon+{\rm gap}(-\delta)\leq \stheta_*$. We first need to bound the quantity 
$\widehat w_N(t;\delta;\epsilon)$ which is tantamount bounding the quantities $\gap(-\delta)$, $\diam(X^{*,\epsilon_1})$ and $\gamma^{(\alpha)}(X^{*,\epsilon_1})$. In the following, $C>0$ is an absolute constant and $\mathfrak{C}>0$ is a constant depending on $\eta_*$, $\vartheta_*$, $\sigma$, $\sigma_*$, $\alpha$, $\kappa$, $\mathfrak{c}$ and $\diam(X^{*})$ that may change from line to line. We use the abbreviation $\phi:=\phi_N(t)$.

From Assumption \ref{assump:LRSS}, \eqref{eq:boundaryptdelta}, \eqref{eq:boundaryptdelta:gap} and $\delta\le\eta_*/2$,
\begin{align}\label{cor:LRSS:random:eq1}
\gap(-\delta)\le \mathfrak{C}\delta^\alpha.
\end{align}    
By  Assumption \ref{assump:LRSS} and $\epsilon_1=\epsilon+\gap(-\delta)\le\vartheta_*$, 
\begin{align*}
\diam(X^{*,\epsilon_1})\le \mathfrak{c}(\epsilon+\gap(-\delta))^\kappa+\diam(X^*)\le
\mathfrak{c}(\epsilon^\kappa+\mathfrak{C}^\kappa\delta^{\alpha\kappa})+\diam(X^*).
\end{align*}
From definition of $\mathfrak{C}_\alpha$, $\phi_N(t)$ and $\widehat w_N(t;\delta;\epsilon)$,
\begin{align}
\widehat{w}_N(t;\delta;\epsilon)&\le C\sigma \diam^\alpha(X^{*,\epsilon_1})\phi
\le \mathfrak{C}\left(\epsilon^{\kappa\alpha}+\delta^{\alpha^2\kappa}+\diam^\alpha(X^*)\right) \phi.\label{cor:LRSS:random:eq2}
\end{align}

\emph{Upper bound on $\check \delta(t;\epsilon)$}. From \eqref{cor:LRSS:random:eq2}, the last inequality defining $S_{N,\eta_*}(t;\epsilon)$ is satisfied if 
$
\delta>\mathsf{A}_N\delta^{\alpha^2\kappa}+\mathsf{B}_N(\epsilon),
$
with 
$
\mathsf{A}_N:=\mathfrak{C}\phi 
$
and
$
\mathsf{B}_N(\epsilon):=\mathfrak{C}(1+\diam^\alpha(X^*)+\epsilon^{\kappa\alpha})\phi.
$
Let
$$
\delta_N(\epsilon):=\max\{2\mathsf{B}_N(\epsilon),(2\mathsf{A}_N)^{1/(1-\alpha^2\kappa)}\}.
$$ 
First, $\delta_N(\epsilon)$ belongs to the set
$
\{\delta:\delta>\mathsf{A}_N\delta^{\alpha^2\kappa}+\mathsf{B}_N(\epsilon)\}
$. 
Second, one has $0<\delta_N(\epsilon)\le\eta_*/2$ if $\mathfrak{C}\phi\le1$ for large enough $\mathfrak{C}>1$. Thirdly, from \eqref{cor:LRSS:random:eq1} and $\epsilon\le\vartheta_*/2$, the constraint $\epsilon+\gap(-\delta_N(\epsilon))\le\vartheta_*$ is satisfied asking for $\mathfrak{C}\phi_N(t)\le1$ for possibly larger $\mathfrak{C}$.  We thus conclude that $\delta_N(\epsilon)\in S_{N,\eta_*}(t;\epsilon)$ implying that
$
    \check\delta(t;\epsilon)=\inf S_{N,\eta_*}(t;\epsilon)\le \delta_N(\epsilon).
$

\emph{Upper bound on $\check r(t;\epsilon_0)$}. Fix $0\le\epsilon_0<\min\{\vartheta_*/4,1/2\}$. Since $\check\delta(t;\epsilon)\le\delta_N(\epsilon)$ and the fact that $\delta\mapsto\widehat w_N(t;\delta;\epsilon)$ is nondecreasing, we get from \eqref{cor:LRSS:random:eq2} and the facts that $\delta_N(\epsilon)\le\eta_*/2$ and $\epsilon+\gap(-\delta_N(\epsilon))\le\vartheta_*$,
\begin{align}
\check w(t;\epsilon)\le\widehat w_N(t;\delta_N(\epsilon);\epsilon)
\le\mathsf{C}_N + \mathsf{D}_N \epsilon^{\kappa\alpha}      
+ \mathsf{E}_N \epsilon^{\kappa^2\alpha^3},\label{cor:LRSS:random:eq3}
\end{align}
with the definitions
$$
\mathsf{C}_N:=\mathfrak{C}\left(1+ \phi^{\alpha^2\kappa}+\mathsf{A}_N^{\alpha^2\kappa/(1-\alpha^2\kappa)}\right)\phi,
$$
as well as
$
\mathsf{D}_N:=\mathfrak{C}\phi
$
and
$
\mathsf{E}_N:=\mathfrak{C}\phi^{1+\alpha^2\kappa}.
$
Moreover, from $\check\delta(t;\epsilon)\le\delta_N(\epsilon)$, $\delta_N(\epsilon)\le\eta_*/2$, \eqref{cor:LRSS:random:eq1} and the fact that  $\delta\mapsto\gap(-\delta)$ is nondecreasing, one has 
\begin{align}
\gap(-\check\delta(t;\epsilon))\le \mathfrak{C}\delta_N^\alpha(\epsilon)
\le \mathsf{F}_N+\mathsf{G}_N\epsilon^{\alpha^2\kappa},\label{cor:LRSS:random:eq4}
\end{align}
with the definitions 
$
\mathsf{F}_N:=\mathfrak{C}\sigma_*^\alpha\phi^\alpha+\mathfrak{C}\mathsf{A}_N^{\alpha/(1-\alpha^2\kappa)}
$
and 
$
\mathsf{G}_N:=\mathfrak{C}\phi^\alpha.
$
From \eqref{cor:LRSS:random:eq3}-\eqref{cor:LRSS:random:eq4}, in upper bounding $\inf R_{N,\eta_*}(t;\epsilon_0)$ it is enough to take the infimum over $\epsilon$ belonging to the set
\begin{align*}
    R:=\left\{\epsilon:
    \epsilon_0+\mathsf{F}_N+\mathsf{G}_N\epsilon^{\alpha^2\kappa}
+2(\mathsf{C}_N+\mathsf{D}_N\epsilon^{\kappa\alpha}
+\mathsf{E}_N\epsilon^{\kappa^2\alpha^3})<\epsilon\le\min\{\frac{\vartheta_*}{2},1\}
    \right\}.
\end{align*}
Define
$$
\epsilon_N:=\max\left\{
2\epsilon_0,8\mathsf{F_N}+16\mathsf{C_N},
(16\mathsf{D}_N)^{1/(1-\alpha\kappa)},
(8\mathsf{G}_N)^{1/(1-\alpha^2\kappa)},
(16\mathsf{E}_N)^{1/(1-\alpha^3\kappa^2)}
\right\}.
$$
It is straightforward to check that one gets $\epsilon_N\in R$ as long as $\epsilon_N<\min\{\vartheta_*/2,1\}$. This requirement follows from   $\epsilon_0<\min\{\vartheta_*/4,1/2\}$ and the fact that $\mathfrak{C}\phi\le1$ for enough large $\mathfrak{C}$. We conclude that
$
    \check r(t;\epsilon_0)=\inf R_{N,\eta_*}(t;\epsilon_0)\le \epsilon_N.
$

\emph{Upper bound on $\check\delta(t)$.} Letting $\epsilon\searrow\check r(t;\epsilon_0)$ and using monotonicity, we obtain from $\check\delta(t;\epsilon)\le\delta_N(\epsilon)$ that $\check\delta(t)\le\delta_N(\check r(t;\epsilon_0))\le\delta_N(\epsilon_N)$. Examining the expressions of $\epsilon_N$ and $\delta_N(\epsilon_N)$ one may check that $\epsilon_N\le\mathfrak{R}_N$ and  $\delta_N(\epsilon_N)\le\mathfrak{C}\phi+\mathfrak{C}\phi\max\{
\epsilon_0^{k\alpha},\phi^{\alpha^2\kappa}\}=:\mathfrak{D}_N$ by enlarging $\mathfrak{C}$ if necessary.

To finalize, Theorem \ref{thm:convexrandomset} and $\check r(t;\epsilon_0)\le\mathfrak{R}_N$ and $\check\delta(t)\le\mathfrak{D}_N$ entail the claim.\qed
\end{proof}

In case of Lipschitz continuity ($\alpha=1$), one has the optimality slackness $\mathfrak{R}_N$ of the order $\epsilon_0+\mathfrak{C}\sqrt{(\mathfrak{C}_\alpha+\log|\mathcal{I}|+t)/N}$ and the feasibility slackness $\mathfrak{D}_N$ of the order $\mathfrak{C}\sqrt{(\mathfrak{C}_\alpha+\log|\mathcal{I}|+t)/N}$. These ``localized rates'' are independent of $\diam(X)$ allowing for an unbounded $X$. They do depend however on the diameter and complexity of $X^{*,\vartheta_*}$. Note that for $\kappa=1/2$ these rates are worse than the case of a fixed feasible set (Proposition \ref{prop:LRSS}). This rate deterioration implied by random constraints is unavoidable in general (see Remark \ref{rem:onedimsimple}).

\section{Application to metric projection problems}\label{sec:simpleresults}

In Section \ref{ss:particular:case} we presented in Propositions \ref{prop:LRSS}-\ref{prop:LRSS:random:set} an application of the localization technique (Theorem \ref{thm:convexrandomset}) in case the solution set satisfies Assumption \ref{assump:LRSS}. Still, in both of these applications, the H\"older modulus is assumed ``uniform'' in the sense that $\sigma^2(\vartheta,\delta)\equiv\sigma^2$ is constant. In this section we present another application where it is important to consider that the H\"older modulus $\sigma^2(\vartheta,\delta)$ varies across the feasible set. The road map will be similar to the proof of Proposition \ref{prop:LRSS:random:set}  with some additional technicalities.

Specifically, we sketch the application of Theorems \ref{thm:generalrandomset}-\ref{thm:convexrandomset} to a simple problem illustrating the difference between the two results. Specifically, we consider a metric projection problem where $X\subset Y\subset \re^d$ and:
\[f_0(x):=\|x-x_0\|^2,\mbox{ with $\|\cdot\|$ the Euclidean norm.}\]
A minimizer $x^*$ of $f_0$ over $X$ corresponds to the metric projection of $x_0$ over $X$. We set $f^*:=R^2$ to be the value of the problem. 

As usual, we assume that $f_0(x) = \espn\,F_0(x,\cdot) = \int_{\Xi}F_0(x,\xi)\,\probn(d\xi)$. Potential examples include:
\begin{enumerate}
\item  \(\Xi = \re_+\times \re\), \(\xi=(\xi[1],\xi[2])^T\sim\probn\) is a random vector with mean $(1,0)^T$, and  \(F_0(x,\xi):= \xi[1]\,\|x-x_0\|^2 + \xi[2];\)
\item \(\Xi = \re^d,\, \xi\sim\probn\) is an isotropic random vector, that is,  satisfying $\esp\langle\xi,x\rangle^2=\Vert x\Vert^2$ for all $x\in\re^d$. In our setting, \(F_0(x,\xi):=\langle \xi,x-x_0\rangle^2\).
\end{enumerate}
In both examples, mild moment conditions on $\probn$ imply that the $\parametersgoodness$-goodness assumption is satisfied with $\alpha=1$ over any bounded set $Z\subset Y$, with a value $\sigma^2=\sigma^2(Z)$ that will in general depend on the set $Z$. Recalling that $\mathbb{B}$ is the unit ball in $\re^d$, we will assume the following condition:
\begin{equation}\label{eq:lipschitzcondition}
\mbox{$\forall r>0$, 
$\parametersgoodness$-goodness holds over $Z=x_0 + r\mathbb{B}$ with } \sigma = c_0\,r,
\end{equation}
where $c_0>0$ is a constant. This condition is compatible with the quadratic growth of $f_0$ in our two examples, with a $c_0$ that depends on $\probn$. For convenience, we assume $\mathcal{I}\neq \emptyset$. 

\begin{remark}The constant $\sigma^2$ in the second example above will inevitably depend on $\|\xi\|^2$. The expectation of $\|\xi\|^2$ is $d$ under our assumptions, which implies $c_0\geq d$ in our assumptions. In this specific setting of $F_0(x,\xi) = \langle x-x_0,\xi \rangle^2$, it has been noticed by Mendelson and others \cite{mendelson2015,mendelson2017} that one can bound the quadratic form $\widehat{F}_0(x)$ from below using very weak assumptions that lead to improved bounds. We will return to this issue in the companion paper \cite{2020oliveira:thompsonII}.
\end{remark}

Before continuing, we note the following direct consequence of strong convexity. Note that: 

\[\forall x\in X\,:\, f(x) = f^* + 2\langle x-x^*,x^*-x_0\rangle + \|x-x^*\|^2\mbox{ and }\langle x-x^*,x^*-x_0\rangle\geq 0.\]
From strong-convexity, 
\begin{align*}
\forall \stheta>0\,:\,\forall x\in X^{*,\stheta},\, \|x-x^*\|\leq  \sqrt{\stheta},
\end{align*}
and, in particular, for some universal constant $C>0$,
\begin{align}
\forall \stheta>0\,:\,\diam(X^{*,\stheta})\leq 2\sqrt{\stheta}\quad\mbox{and}\quad \gamma_2^{(1)}(X^{*,\stheta})\leq C\sqrt{d\stheta},
\label{eq:approx:sol:diam:gamma:func}
\end{align}
using Dudley's bound in \eqref{eq:simpleboundgamma2}. We note that the second bound might be far from sharp in several examples.

\subsection{Application of Theorem \ref{thm:generalrandomset}} We make the $\parametersgoodness$-goodness assumption over $Z=Y$ (cf. Assumption \ref{assump:goodFi}). We treat $\sigma_*$, the parameter $c_0$ in (\ref{eq:lipschitzcondition}) as universal constants, and use $C$ to denote a  universal constant that might change from line to line. We set $\epsilon_0=0$, as we are interested in exact minimizers of the SAA problem.

Recalling that $R^2:=\inf_{x\in X}\|x-x_0\|^2$, we obtain that $\sigma\leq C\,(R+{\rm diam}(Y))$ (cf. (\ref{eq:lipschitzcondition})). The parameters 
$\widehat{r}_N(t)$ and $\widehat{\delta}_N(t)$ satisfy:
\[\widehat{r}_N(t)\leq \widehat{\delta}_N(t)\le\frac{C\,(R+{\rm diam}(Y))\,(\diam(Y)\sqrt{t + \log(|\mathcal{I}|+1)} + \gamma_2^{(1)}(Y))}{\sqrt{N}}.\]

If $x^*$ belongs to the relative interior of $X$ in $Y$ (ie. there exists $r>0$ with $(x^*+r\mathbb{B})\cap Y\subset X$, then
\[{\rm gap}(-\widehat{\delta}_N(t))\leq \frac{C\,(R+{\rm diam}(Y))^2\,(\diam(Y)^2(t + \log(|\mathcal{I}|+1)+ (\gamma_2^{(1)}(Y))^2)}{N}.\]
In general, using (\ref{eq:holdercont}), the metric regularity condition (Assumption \ref{assump:MRF}) and the above estimate on 
$\sigma$ to deduce that\footnote{Indeed, if $x_{\delta}\in X_{\delta}$ is the metric projection of $x^*_\delta$ onto $X$ for some $\delta>0$, by (\ref{eq:holdercont}) and Assumption \ref{assump:MRF}, we have $f^*-f_{\delta}^*\le f(x_\delta)-f(x_\delta^*)\le\sigma\dist(x^*_\delta,X)\le\sigma\mathfrak{c}\delta$.}
\begin{align*}
{\rm gap}(\widehat{\delta}_N(t))&\leq  \frac{C\,(R+{\rm diam}(Y))^2\,(\diam(Y)\sqrt{t + \log(|\mathcal{I}|+1)} + \gamma_2^{(1)}(Y))}{\sqrt{N}},
\end{align*}
where now $C$ depends on $\mathfrak{c}$ from Assumption \ref{assump:MRF} as well. In order to bound 
${\rm gap}(-\widehat{\delta}_N(t))$, we assume a slightly stronger version of Assumption \ref{assump:MRF}:
\begin{equation*}
\exists\eta_*>0,\forall 0\le\delta\le\eta_*,\forall x\in Y,\quad
\dist(x,X_{-\delta})\le \mathfrak{c}\sup_{i\in\mathcal{I}}[f_i(x)+\delta]_+.
\end{equation*}
Similarly,
\begin{align*}
{\rm gap}(-\widehat{\delta}_N(t))&\leq  \frac{C\,(R+{\rm diam}(Y))^2\,(\diam(Y)\sqrt{t + \log(|\mathcal{I}|+1)} + \gamma_2^{(1)}(Y))}{\sqrt{N}},
\end{align*}
where now $C$ depends on $\mathfrak{c}$ and $\eta_*$. By Theorem \ref{thm:generalrandomset}, we obtain that with probability $\geq 1-e^{-t}-(2|\mathcal{I}+1|)\rho$, 
\begin{eqnarray*}\|\widehat{x}^*-x_0\|^2 &\leq & R^2  +\frac{C\,(R+{\rm diam}(Y))^2\,(\diam(Y)\sqrt{t + \log(|\mathcal{I}|+1)} + \gamma_2^{(1)}(Y))}{\sqrt{N}};\\
\max_{i\in\mathcal{I}}f_i(\widehat{x}^*)&\leq & \frac{C\,(R+{\rm diam}(Y))^2\,(\diam(Y)\sqrt{t + \log(|\mathcal{I}|+1)} + \gamma_2^{(1)}(Y))}{\sqrt{N}}.
\end{eqnarray*}
The bounds above are of the order $N^{-1/2}$, which coincides with what comes from asymptotic analyses. Other interesting aspects of our results are the explicit dependence on ${\rm diam}(Y)$, $\gamma_2^{(1)}(Y)$ and $R$. In the next subsection, we show that these bounds can be refined significantly under the assumptions of Theorem \ref{thm:convexrandomset}.

\subsection{Application of Theorem \ref{thm:convexrandomset}}\label{sub:applyconvexrandomset}

We now work under the assumptions of Theorem \ref{thm:convexrandomset} combined with our discussion in the beginning of the section. We treat $\sigma_*$, $\eta_*$, $\stheta_*$ and the constant $c_0$ in (\ref{eq:lipschitzcondition}) as absolute constants, and use $C,C_0>0$ to denote generic constants depending only on $\sigma_*,\eta_*,\stheta_*$, 
$\diam(X^{*,\vartheta_*})$ and $c_0$. In particular, their precise values may be different in each occurrence. We assume without loss on generality that $R\ge1$. We set 
$\epsilon_0=0$, as we are interested in exact minimizers of the SAA problem. Fix also $t>0$.

Let $0<\epsilon\leq \min\{\stheta_*/2,1\}$ and $0<\delta\le\min\{\eta_*/2,1\}$ with $\epsilon+{\rm gap}(-\delta)\leq \stheta_*$. For simplicity let $\epsilon_1:=\epsilon+\gap(-\delta)$.  We first need to bound the quantity 
$\widehat w_N(t;\delta;\epsilon)$ which is tantamount bounding the quantities
$\diam(X^{*,\epsilon_1})$, $\gamma^{(1)}(X^{*,\epsilon_1})$ and 
$\sigma(\epsilon_1;\delta)$. 

\begin{description}
\item[\emph{Bound on $\gap(-\delta)$}:]
If $x_{-\delta}$ is the metric projection of $x^*$ onto $X_{-\delta}^{*,\stheta_*}$, \eqref{eq:boundaryptdelta} above implies $\Vert x_{-\delta}-x^*\Vert\le C\delta$, with $C>0$ depending on $\diam(X^{*,\vartheta_*})/\eta_*$. Hence, as $\Vert x_0-x^*\Vert=R$, we have that $x^*,x_{-\delta} $ lies in $x_0+r\mathbb{B}$ with $r\le C(R+\delta)$. From the goodness assumption in \eqref{eq:lipschitzcondition} with a Lipschitz modulus 
$\sigma =c_0 r$ over $x_0+r\mathbb{B}$, we get 
\begin{align}
\gap(-\delta)\le f(x_{-\delta})-f^*\le C(R+\delta)\Vert x_{-\delta}-x^*\Vert\le C(R\delta+\delta^2).\label{eq:gap:delta}
\end{align}    

\item[\emph{Bound on $\diam(X^{*,\epsilon_1})$}:]
By Lemma \ref{lem:robinson2}, \eqref{eq:approx:sol:diam:gamma:func} and the previous bound on $\gap(-\delta)$, we have 
\begin{align*}
\diam(X_\delta^{*,\epsilon+\gap(-\delta)})&\le2\diam(X^{*,\epsilon+\gap(-\delta)})\\
&\le C\sqrt{\epsilon+\gap(-\delta)}\le C(\sqrt{\epsilon}+\sqrt{R\delta}+\delta).
\end{align*}

\item[\emph{Bound on $\sigma(\epsilon_1;\delta)$}:]
Recall that $\sigma(\epsilon+\gap(-\delta);\delta)$ is the Lipschitz constant over the set $X_\delta^{\epsilon+\gap(-\delta)}$. Of course, $x^*\in X_\delta^{\epsilon+\gap(-\delta)}$ and we already shown $\diam(X_\delta^{*,\epsilon+\gap(-\delta)})\le C(\sqrt{\epsilon}+R+\delta)$. Hence, $\Vert x_0-x^*\Vert=R$ and triangle inequality yield $X_\delta^{\epsilon+\gap(-\delta)}\subset x_0+r\mathbb{B}$ with $r=C(\sqrt{\epsilon}+R+\delta)$. The goodness assumption in \eqref{eq:lipschitzcondition} thus implies 
$$
\sigma(\epsilon+\gap(-\delta);\delta)\le C(\sqrt{\epsilon}+R+\delta).
$$
\end{description}

Now, let us define 
$$
\phi:=\sqrt{\frac{d+\log|\mathcal{I}|+t}{N}}.
$$
Using the above bounds, Dudley's bound \eqref{eq:simpleboundgamma2}, which yields
$$
\gamma^{(1)}(X^{*,\epsilon+\gap(-\delta)})\le C\sqrt{d}(\sqrt{\epsilon}+\sqrt{R\delta}+\delta),
$$
and a simple but tedious computation\footnote{Using 
$\sqrt{\epsilon R\delta}\le2\epsilon+2R\delta$ and $\epsilon,\delta\le1$.}, we obtain
\begin{align}
\widehat{w}_N(t;\delta;\epsilon)&\le C\phi\sigma(\epsilon_1;\delta)\diam(X^{*,\epsilon_1})\le C(\sqrt{\epsilon}+R+\delta)(\sqrt{\epsilon}+\sqrt{R\delta}+\delta)\nonumber\\
&\le C\phi(R^{3/2}\sqrt{\delta}+\delta^{2})+C\phi(R\sqrt{\epsilon}+\epsilon),\label{eq:wN}
\end{align}
for all $0<\epsilon\leq \min\{\stheta_*/2,1\}$ and $0<\delta\le\min\{\eta_*/2,1\}$ with $\epsilon+{\rm gap}(-\delta)\leq \stheta_*$.

\emph{Bound on $\check\delta(t;\epsilon)$.} With a bound on $\widehat{w}_N(t;\delta;\epsilon)$, we may obtain a sufficient upper bound on $\check\delta(t;\epsilon)$ using the definition of $S_{N,\eta_*}(t;\epsilon)$. Using that $\delta^2\le1$, for the constraint $$
\delta>\widehat{w}_N(t;\delta;\epsilon)+\sigma_*\sqrt{6\frac{1+\log(2(|\mathcal{I}|+1))+t}{N}},
$$
to hold, it suffices 
\begin{align}
\delta>\mathsf{A}\sqrt{\delta}+\mathsf{B}(\epsilon),\label{eq:abc}
\end{align}
with 
$
\mathsf{A}:=C\phi R^{3/2}
$
and
$
\mathsf{B}(\epsilon):=C\phi(R\sqrt{\epsilon}+\epsilon+1).
$
First, for \eqref{eq:abc} to hold, it thus suffices to choose 
$$
\delta(\epsilon):=\max\left\{4\mathsf{A}^2,2\mathsf{B}(\epsilon)\right\}.
$$
Second, one has $0<\delta(\epsilon)\le\eta_*/2$ if $C_0 R^{3/2}\phi\le1$ for large enough $C_0>1$. Thirdly, from \eqref{eq:gap:delta} and $\epsilon\le\vartheta_*/2$, the constraint $\epsilon+\gap(-\delta_N(\epsilon))\le\vartheta_*$ is satisfied asking for possibly larger $C_0$. We thus conclude that $\delta(\epsilon)\in S_{N,\eta_*}(t;\epsilon)$ implying that
\begin{align}
    \check\delta(t;\epsilon)=\inf S_{N,\eta_*}(t;\epsilon)\le \delta(\epsilon)\le C\phi^2 R^3+C\phi(R\sqrt{\epsilon}+\epsilon+1),\label{eq:check:delta}
\end{align}
for all $0<\epsilon\le\min\{\vartheta_*/2,1\}$.

\emph{Bound on $\check r(t;0)$.} Let for all  $0<\epsilon\le\min\{\stheta_*/2,1\}$. Recall $\check w(t;\epsilon)=\widehat w_N(t;\check\delta(t;\epsilon);\epsilon)$. We now pursue an upper bound on $\check r(t;0)$ by checking the definition of $R_{N,\eta_*}(t;0)$. After some computations, using that $0<\delta(\epsilon)\le\min\{\eta_*/2,1\}$ and $\epsilon+{\rm gap}(-\delta(\epsilon))\leq \stheta_*$ it follows from monotonicity and 
\eqref{eq:check:delta} and \eqref{eq:wN} that 
\begin{align*}
\check w(t;\epsilon)
&\le C\phi^{3/2}R^2\epsilon^{1/4}+C(\phi^{3/2}R^{3/2}+\phi R)\sqrt{\epsilon}
+\frac{\epsilon}{4}\\
&+C(\phi^3+\phi^5R^6+\phi^{3/2}R^{3/2}+\phi^2R^3),
\end{align*} 
where we used that $0<\epsilon\le1$ and  $C(\phi^3R^2+\phi)\epsilon+C\phi^3\epsilon^2\le\epsilon/4$ by enlarging $C_0$ if necessary. Moreover, by \eqref{eq:gap:delta} and \eqref{eq:check:delta} we get
\begin{align*}
\gap(-\check\delta(t;\epsilon))
&\le C\phi R^2\sqrt{\epsilon}+\frac{\epsilon}{4}
+C(\phi^4R^6+\phi^2R^4+R\phi+\phi^2),
\end{align*}
using $C(R\phi+\phi^2R^2)\epsilon+C\phi^2\epsilon^2\le\epsilon/4$ for large enough $C_0$. From the two previous displays, in order to have 
\begin{align}
2\check w(t;\epsilon)+\gap(-\check\delta(t;\epsilon))<\epsilon,\label{eq:check:w:gap:epsilon}
\end{align}
it is enough that
\begin{align*}
C\phi^{3/2}R^2\epsilon^{1/4}&<\epsilon/10,\\
C(\phi^{3/2}R^{3/2}+\phi R+\phi R^2)\sqrt{\epsilon}&<\epsilon/10,\\
C(\phi^4R^6+\phi^5R^6+\phi^2R^4+\phi^2R^3+\phi^{3/2}R^{3/2}+R\phi
+\phi^3+\phi^2)&<\epsilon/10.
\end{align*}
For the above conditions to hold, one may check that it is enough to have
\begin{align}
\epsilon_N:=C[\phi^{3/2}R^{3/2}+(R^4+1)\phi^2].\label{eq:epsilon:metric}
\end{align}
By enlarging $C_0$ if necessary we may also guarantee the additional constraint $0<\epsilon_N\le\stheta_*/2$ as required in $R_{N,\eta_*}(t;0)$. In conclusion, $\epsilon_N\in R_{N,\eta_*}(t;0)$ and hence
\begin{align*}
\check r(t;0)=\inf R_{N,\eta_*}(t;0)\le \epsilon_N
\le C[\phi^{3/2}R^{3/2}+(R^4+1)\phi^2].
\end{align*}

\paragraph{Bound on $\check\delta(t)$.} Setting $\epsilon\searrow\check r(t;0)$ in \eqref{eq:check:delta} and using \eqref{eq:epsilon:metric},
simple calculation yields
\begin{align*}
\check\delta(t)&\le C\phi^{7/4}R^{7/4}+C\phi^{5/2}R^{3/2}
+C\phi^3R^4+C\phi^2 R+C\phi^2R^3+C(\phi+\phi^3)\\
&\le C(\phi^2R^3+\phi^2R+\phi^{1.75}R^{1.75}+\phi),
\end{align*} 
where we used that $C_0\phi R^{3/2}\le1$ for large enough $C_0\ge1$.

Recall $\phi=\sqrt{\frac{d+\log|\mathcal{I}|+t}{N}}$. From Theorem \ref{thm:convexrandomset} and the fact that $\gap(-\check\delta(t))\le\check r(t;0)$ by  \eqref{eq:check:w:gap:epsilon}, we conclude that, for $N$ large enough so that $C_0\phi R^{3/2}\le1$, with probability $\geq 1-e^{-t}-(2|\mathcal{I}|+2)\rho$, 
\begin{eqnarray*}\|\widehat{x}^*-x_0\|^2 &\leq & R^2+ C[\phi^{3/2}R^{3/2}+(R^4+1)\phi^2];\\
\dist(\widehat{x}^*,X) &\leq & C\phi[1+\phi(R+R^3)]+C\phi^{1.75}R^{1.75}.
\end{eqnarray*}
For large $d,N$, $f(\widehat{x}^*)-f(x^*)$ decays like $R^{3/2}(d/N)^{3/2}+(R^4+1)d/N$. Note that it depends on $\diam(X^{*,\vartheta_*})$ but not on the diameters of $X$ nor $Y$.

Now assume additionally that 
\begin{equation}\label{eq:interiorpoint}
x_0=x^*\in X_{-\delta_*}\mbox{ with }\delta_* = C_*\sqrt{\frac{1+t+\log(1+|\mathcal{I}|)}{N}}
\end{equation} 
with sufficiently large $C_*$; ie. $x_0$ is ``sufficiently interior"~to $X$. Then we have $R=0$ and ${\rm gap}(-\delta)=0$ for $0\leq \delta\leq \delta_*$. One can see that, in this case, the dependence on $0\leq \delta\leq\delta_* $ disappers in the bounds related to ${\rm gap}(-\delta)$. Some calculations then improve our high-probability bound on $\|\widehat{x}^*-x_0\|^2$ to:
\begin{eqnarray}\label{eq:improvedrate} 
\|\widehat{x}^*-x_0\|^2&\leq & C\left(\frac{d+t+\log(|\mathcal{I}|+1)}{N}\right).
\end{eqnarray}
whenever $N\geq C_0(d+t+\log(|\mathcal{I}|+1))$. This is the kind of fast rate expected in strongly convex problem. However, such an improvement requires that $x_0$ be a ``sufficiently interior" point of $X$; see Remark \ref{rem:onedimsimple} below.  

Finally, when $N<C_0(d+t+\log(|\mathcal{I}|+1))$ and (\ref{eq:interiorpoint}) does not hold, we may still obtain a bound by using:
\[\gamma_2^{(1)}(X^{*,\epsilon + {\rm gap}(-\delta)})\leq \gamma_2(X^{*,\stheta_*})\]
for small enough $\epsilon,\delta$. This leads to nontrivial bounds whenever $\gamma_2(X^{*,\stheta_*})\ll \sqrt{d}$ (eg. if $X^{*,\stheta_*}$ is contained in a small simplex).  

\begin{remark}\label{rem:onedimsimple} We observe that even in one dimension we may expect fluctuations of order $R^2/\sqrt{N}$ on $\|\widehat{x}^*-x_0\|^2-R^2$ in metric projection problems, when $x_0$ lies outside the feasible set. Assume $\Xi=Y=\re$, $\xi\sim \probn$ is exponential with parameter $1$ (that is, $\prob\{\xi>t\} =e^{-t}$ for all $t>0$), $F_0(x,\xi)=f_0(x)=x^2$, $\mathcal{I}=\{1\}$ and:
\[F_1(x,\xi) = R - \xi\,x.\] 
The solutions to the ideal problem and SAA are $x^*=R$ and $\widehat{x}^* = R/\bar{\xi}_N$,  
where $\bar{\xi}_N$ is the sample average of the $\xi_k$. Using the Central Limit Theorem for $\bar{\xi}_N$, one can show that  $\sqrt{N}((\widehat{x}^*)^2 - R^2)$ has a Gaussian limit with standard deviation $2R^2$ when $N\to +\infty$.\end{remark}

\section{Concentration inequalities for sample averages}\label{sec:concentration}

We present a novel concentration inequality. The proof of Theorems \ref{thm:generalrandomset} and \ref{thm:convexrandomset}
rely on this tool which may be of independent interest in stochastic optimization.

\begin{theorem}[Proof below]\label{thm:concentration} Suppose $(\mathcal{M},\dist)$ is a totally bounded metric space. Assume \[(\Xi,\sigma(\Xi)),\,\probn,\,\{\xi_k\}_{k=1}^N\mbox{ and }\xi\sim \probn\] are as in Section \ref{sec:basic} and $G:\mathcal{M}\times \Xi\to \re$ is a measurable function with $\espn\,|G(x_0,\cdot)|<+\infty$ for some $x_0\in \mathcal{M}$. Assume additionally that there exists a measurable function $\mathsf{L}:\Xi\to \re_+$ with $\espn\,\mathsf{L}^2\leq \nu^2<+\infty$ and a constant $0<\alpha\leq 1$ such that:  
\[\mbox{ for $\probn$-a.e. $\xi\in\Xi$, }\forall x,x'\in\mathcal{M}\,:\,|G(x,\xi) - G(x',\xi)|\leq \mathsf{L}(\xi)\,\dist(x,x')^{\alpha}.\]
Write:\[\Delta G := \sup_{x\in \mathcal{M}}|(\widehat{\espn} - \espn)(G(x,\cdot) - G(x_0,\cdot))|\]
and assume:
\[\prob\{\widehat{\espn}\mathsf{L}^2(\cdot)>2\nu^2\}\leq \rho\in [0,1].\]
Then, for any $t\geq 0$:
\[\prob\left\{\widehat{\espn}\mathsf{L}(\cdot)^2\leq 2\nu^2,\,\Delta G >\nu\frac{ \cgammaconc\,\gamma^{(\alpha)}_2(\mathcal{M},\dist) +\cdiamconc\diam(\mathcal{M})^{\alpha}\sqrt{1+t}  }{\sqrt{N}}\right\}\leq e^{-t}.\]
\end{theorem}

Notice that, if $N$ grows, $\widehat{\espn}\mathsf{L}^2\to {\espn}\mathsf{L}^2\leq \nu^2$ almost surely. Therefore, we expect the probability of $\widehat{\espn}\mathsf{L}(\cdot)^2\leq 2\nu^2$ to be large when $N$ is large. The above theorem shows that on the event that $\widehat{\espn}\mathsf{L}(\cdot)^2\leq 2\nu^2$, the likelihood of $\Delta G$ being large is exponentially small. 

To prove this result, we will use the next lemma. It is a simple consequence of a much more general result of Panchenko \cite{panchenko}. 

\begin{lemma}[Proof in Appendix]\label{lem:panchenko} Assume $Z_1,\dots,Z_N$ are i.i.d. random variables with finite second moments. Then
\[\prob\left\{\frac{1}{N}\sum_{i=1}^N(Z_i - \esp[Z_1])>\sqrt{\frac{2(1+t)}{N}\left(\var[Z_1] + \frac{1}{N}\sum_{i=1}^N(Z_i-\esp[Z_1])^2\right)}\right\}\leq e^{-t}.\]
In particular, if $\var[Z_1]\leq \nu^2$, 
\[\prob\left\{\frac{1}{N}\sum_{i=1}^N(Z_i - \esp[Z_1])>\nu\sqrt{\frac{6(1+t)}{N}},\;\frac{1}{N}\sum_{i=1}^N(Z_i-\esp[Z_1])^2\leq 2\nu^2\right\}\leq e^{-t}.\]
\end{lemma}

This remarkable inequality by Panchenko shows that averages of $Z_i$, when normalized by an {\em empirical term}, have sub-Gaussian tails under extremely weak assumptions. We will use this both to prove Theorem \ref{thm:concentration} and to control fluctuations of other random variables in the proofs of Theorems \ref{thm:generalrandomset} to \ref{thm:convexrandomset}.

\begin{proof}[of Theorem \ref{thm:concentration}] In this proof we use a combination of generic chaining (as encapsulated by Theorem \ref{theorem:gc}) and Panchenko's self-normalized concentration inequality (Lemma \ref{lem:panchenko}). 

We begin by noting that, since $\espn\mathsf{L}^2(\cdot)\leq \nu^2$, for any $u>0$, 
\begin{align*}
\prob\left\{\widehat{\espn}\,\mathsf{L}^2(\cdot)\le2\nu^2,
\Delta G>u\,\sqrt{\frac{6\nu^2}{N}}\right\}
&\le \prob\left\{
(\widehat{\espn}+\espn)\,\mathsf{L}^2(\cdot)\le3\nu^2,
\Delta G>u\,\sqrt{\frac{6\nu^2}{N}} 
\right\}\\
&\le \prob\left\{
\Delta G>u\,\sqrt{\frac{2(\widehat{\espn}+\espn)\,\mathsf{L}^2(\cdot)}{N}} 
\right\}.
\end{align*}
Letting $t\geq 0$ and
\[u=u_t:=\cgamma\gamma^{(\alpha)}_2(\mathcal{M},\dist) +\cdiam\diam(\mathcal{M})^{\alpha}\sqrt{1+t},\]
we see that it suffices to show that: 
\begin{equation}\label{eq:sufficientconcentration}{\bf Sufficient:}\;\prob\left\{\Delta G>u_t\sqrt{\frac{2(\widehat{\espn} + \espn)\,\mathsf{L}^2(\cdot)}{N}}\right\}\leq e^{-t} .\end{equation}
To prove (\ref{eq:sufficientconcentration}), we will use our ``generic chaining"~ bound, Theorem \ref{theorem:gc}. For each $x\in \mathcal{M}$, define the random quantity:
\[Y_x:= \frac{(\widehat{\espn} - \espn)\,G(x,\cdot)}{\sqrt{\frac{2(\widehat{\espn} + \espn)\,\mathsf{L}^2(\cdot)}{N}}},\]
when the denominator is $\neq 0$, or $Y_x=0$ otherwise.  Note that:
\[\Delta G = \sqrt{\frac{2(\widehat{\espn} + \espn)\,\mathsf{L}^2(\cdot)}{N}}\,\sup_{x\in\mathcal{M}}|Y_x - Y_{x_0}|.\]
If we can show that:
\[{\bf Goal:}\;\forall x,x'\in\mathcal{M},\,\forall t\geq 0\,:\,\prob\{Y_x - Y_{x'}\geq \sqrt{2(1+t)}\dist(x,x')^{\alpha}\}\le e^{-t},\]
then Theorem \ref{theorem:gc} gives us (\ref{eq:sufficientconcentration}). 

To obtain our goal, we fix $x,x'$ and $t$. We will apply Panchenko's inequality (Lemma \ref{lem:panchenko}) to the i.i.d. random variables:
\[Z_k:= G(x,\xi_k) -  G(x',\xi_k)\,(k\in [N]),\]
so that:
\begin{equation}\label{eq:relationZG}(\widehat{\espn} - \espn)\,(G(x,\cdot) - G(x',\cdot)) = \frac{1}{N}\sum_{k=1}^N(Z_k - \esp[Z_k]).\end{equation}
To apply Lemma \ref{lem:panchenko}, we will estimate the terms $\var[Z_1]$ and $(Z_k-\esp[Z_k])^2$ appearing in that bound. Note that:
\begin{eqnarray*}\var[Z_1] + (Z_k-\esp[Z_k])^2 &=& \esp[Z_k^2] + Z_k^2 - 2Z_k\esp[Z_k]\\ &\leq & \esp[Z_k^2] + Z_k^2 + 2\sqrt{Z^2_k\esp[Z^2_k]}\\ \mbox{($2\sqrt{xy}\leq x+y$ for all $x,y\in\re_+$)} &\leq & 2(\esp[Z_k^2] + Z_k^2).\end{eqnarray*}
Now, by our assumptions:
\[|Z_k| = |G(x,\xi_k) -  G(x',\xi_k) |\leq \mathsf{L}(\xi_k)\,\dist(x,x')^\alpha,\]
therefore:
\[ \var[Z_1] + \frac{1}{N}\sum_{k=1}^N(Z_k-\esp[Z_k])^2\leq 2\,[(\widehat{\espn}+\espn)\,\mathsf{L}(\cdot)^2]\,\dist(x,x')^{2\alpha}.\]
We may finally apply Panchenko's inequality and deduce the following bound:
\begin{eqnarray*}\prob\left\{\frac{1}{N}\sum_{k=1}^N(Z_i-\esp[Z_i])\geq 2\sqrt{\frac{(1+t)(\widehat{\espn}+\espn)\,\mathsf{L}(\cdot)^2}{N}}\dist(x,x')^{\alpha}\right\} \leq e^{-t}.\end{eqnarray*}
This implies our goal once we combine it with (\ref{eq:relationZG}) and the definition of $Y_x$.
\qed
\end{proof}

\section{Deviation and localization arguments}\label{sec:deviationlocalization}

This section compiles a series of deterministic results on how the SAA differs from the ideal optimization problem. We consider general sets and functions in \S \ref{sub:deviationgeneral} and the convex case in \S \ref{sub:deviationconvex}. We are particularly careful to distinguish lower and upper tails in our bounds, as lower tails can be much better behaved than upper tails. One such setting will be explored in the companion paper \cite{2020oliveira:thompsonII}.

\begin{remark} The results in this section are purely deterministic in the sense that we do not need \eqref{eq:mean:fi} and \eqref{eq:mean:hat:fi} to hold. We simply need to assume that $Y$ is as given in Section \ref{sec:basic}; that $f_i,\widehat{F}_i:Y\to \re$ are functions (with $i\in\mathcal{I}_0$) and that $X$, $\widehat{X}$, etc are defined in terms of the $f_i$ and $\widehat{F}_i$ as prescribed in Section \ref{sec:basic}. Our results will be the most interesting when the $\widehat{F}_i$ are good approximations to the respective $f_i$. \end{remark}

For convenience, we introduce the following notation. Given $x,y\in Y$, $Z\subset Y$ and $i\in\mathcal{I}_0$:
\begin{eqnarray}\label{eq:defDeltaone}\widehat{\Delta}_i(x)&:=&\widehat{F}_i(x) - f_i(x);\\ 
\label{eq:defdeltatwo}\widehat{\Delta}_i(y;x)&:=& \widehat{\Delta}_i(y)  - \widehat{\Delta}_i(x).\end{eqnarray}

\subsection{General sets and functions}\label{sub:deviationgeneral}

The next lemma is quite straightforward.

\begin{lemma}\label{lem:generaldeviations} Given $\delta,\epsilon,\epsilon_0\geq 0$, assume $|\widehat{\Delta}_i(y)|\leq \delta$ for all $y\in Y$ and $i\in\mathcal{I}$ and also that $|\widehat{\Delta}_0(y;x^*)|\leq \epsilon$ for all $y\in Y$. Then we have the following. 
\begin{enumerate}
\item $X_{-\delta}\subset \widehat{X}\subset X_{\delta}$.
\item $\widehat{X}^{*,\epsilon_0}\subset X_{\delta}^{*,\epsilon_0+2\epsilon+{\rm gap}(-\delta)}$.
\item $|\widehat{F}^* - f^*|\leq \delta +  \max\{2\epsilon+{\rm gap}(-\delta),\gap(\delta)\}$.
\end{enumerate}\end{lemma}
\begin{proof}The first item is an immediate consequence of the fact that $-\delta\leq \widehat{\Delta}_i(y)\leq \delta$ for each $i\in\mathcal{I}$. 

Let $\widehat{x}\in \widehat{X}^{*,\epsilon_0}$ and  $x^*_{-\delta}\in {\rm arg min}_{x\in X_{-\delta}}f(x)$, so that $f(x^*_{-\delta})=f^*_{-\delta} = f^* + {\rm gap}(-\delta)$. By item 1, $x_{-\delta}^*\in \widehat{X}$, so
\[f(\widehat{x}) + \widehat{\Delta}_0(\widehat{x})= \widehat{F}(\widehat{x})\leq \widehat{F}^*+\epsilon_0\leq \widehat{F}(x^*_{-\delta})+\epsilon_0= f^*_{-\delta} + \widehat{\Delta}_0(x^*_{-\delta})+\epsilon_0.\]
Therefore, for any $\widehat{x}\in \widehat{X}^{*,\epsilon_0}\subset X_\delta$, 
\begin{eqnarray}
f(\widehat x) -f^*&\leq f^*_{-\delta}-f^* +  \widehat{\Delta}_0(x^*_{-\delta})+\epsilon_0 - \inf_{y\in X_\delta} \widehat{\Delta}_0(y)\nonumber\\
&\leq {\rm gap}(-\delta)+\epsilon_0 + \sup_{y\in X_{\delta}}|\widehat{\Delta}_0(y;x^*_{-\delta})|\nonumber\\
&\leq {\rm gap}(-\delta)+\epsilon_0 + 2\sup_{y\in X_{\delta}}|\widehat{\Delta}_0(y;x^*)|,\label{lem:generaldeviations:eq1}
\end{eqnarray}
which gives item $2$. 

For item 3, we assume for simplicity that some $\widehat{x}\in \widehat{X}$ achieves the minimum of $\widehat{F}$: $\widehat{F}(\widehat{x})=\widehat{F}^*$. In this case, 
\begin{eqnarray*}\nonumber |\widehat{F}^*-f^*|& =  & |\widehat{F}(\widehat{x}) - f(x^*)| \\ \nonumber &\leq & |\widehat{F}(\widehat{x})-f(\widehat{x})| + |f(\widehat{x})- f(x^*)| \\ &\leq & \sup_{x\in X_{\delta}}|\widehat{\Delta}_0(x)| + |f(\widehat{x})- f(x^*)|.
\end{eqnarray*}
In one hand $f(\widehat x)-f^*$ is upper bounded by \eqref{lem:generaldeviations:eq1} (with $\epsilon_0=0$). Since $\widehat x\in X_\delta$, a lower bound is given by
$
    f^*-f(\widehat x)=\gap(\delta)+f_\delta^*-f(\widehat x)
    \le\gap(\delta),
$
finishing the proof.
\qed
\end{proof}

\subsection{Convex sets and functions}\label{sub:deviationconvex}

We now consider the convex setting with ``localized"~bounds. 

\begin{proposition}\label{lem:convexlocalized}Assume that $Y$ is convex and closed and that the functions $\{f_i\}_{i\in\mathcal{I}_0}$ and $\{\widehat{F}_i\}_{i\in\mathcal{I}_0}$ are all convex and continuous. Given 
${\delta^{\circ}},\delta>0$, assume $x_{-\delta^{\circ}}\in X_{-\delta^{\circ}}$. Fix $\epsilon\geq f(x_{-\delta^{\circ}})-f^*$ and $\epsilon_0>0$. Let 
$\hat\epsilon:=\widehat{F}(x_{-\delta^\circ}) - \widehat{F}^*+\epsilon_0
$. 
If the following three conditions hold:
\begin{align}\label{eq:controlconstraintsinterior} 
\forall i\in\mathcal{I},&&\widehat{\Delta}_i(x_{-\delta^{\circ}})&\leq \delta^{\circ};\\ 
\label{eq:controlconstraintsboundary}
\forall i\in\mathcal{I},&&\inf_{x\in X^{*,\epsilon}_{\delta,{\rm act}(i)}\cap \widehat{X}^{*,\hat\epsilon}} \widehat{\Delta}_i(x)&>  -\delta;\\
\label{eq:controlobjective} &&\inf_{x\in X^{*,=\epsilon}_{\delta}\cap \widehat{X}^{*,\hat\epsilon}} \widehat{\Delta}_0(x;x_{-\delta^{\circ}}) &> -(\epsilon - (f(x_{-\delta^{\circ}}) - f^*)-\epsilon_0),
\end{align}
then:   
\begin{enumerate}
\item \(\widehat{X}^{*,\epsilon_0}\subset X^{*,\epsilon}_\delta,\) or equivalently, any $x\in\widehat{X}$ with $\widehat{F}(x)\leq \widehat{F}^* + \epsilon_0$ satisfies $x\in Y$, $\max_{i\in\mathcal{I}}f_i(x)\leq \delta$ and $f(x)\leq f^*+\epsilon$;
\item The values of the SAA and the ideal problem satisfy \[|\widehat{F}^*- f^*|\leq |\widehat{\Delta}_0(x^*)| + \sup_{x\in X^{*,\epsilon}_{\delta}}|\widehat{\Delta}_0(x;x^*)| + \max\{\epsilon,\gap(\delta)\}.\]
\end{enumerate}\end{proposition}

\begin{proof}The proof consists of three main steps. In the {\em first step}, we show that assumption (\ref{eq:controlconstraintsinterior}) implies $x_{-\delta^\circ}\in \widehat{X}\cap X^{*,\epsilon}$. In the {\em second step}, we show that {if} (\ref{eq:controlconstraintsinterior}) holds and there exists a point $z\in \widehat{X}^{*,\epsilon_0}\backslash X^{*,\epsilon}_\delta$, {\em then} one of (\ref{eq:controlconstraintsboundary}) or (\ref{eq:controlobjective}) {\em cannot} hold. In contrapositive form, the second step implies that, if we assume the three conditions  (\ref{eq:controlconstraintsinterior}), (\ref{eq:controlconstraintsboundary}) and (\ref{eq:controlobjective}), then $\widehat{X}^{*,\epsilon_0}\subset X^{*,\epsilon}_\delta$. Finally, the {\em third step} proves the inequality for $|\widehat{F}^*-f^*|$.

\paragraph{First step.} Assume (\ref{eq:controlconstraintsinterior}). We argue that $x_{-\delta^{\circ}}\in\widehat{X}$. To see this we first observe that $x_{-\delta^{\circ}}\in Y$. Moreover, for all $i\in\mathcal{I}$, $f_i(x_{-\delta^{\circ}})\leq -\delta^{\circ}$, so
\[\widehat{F}_i(x_{-\delta^{\circ}})\leq -\delta^{\circ} + \widehat{\Delta}_i(x_{-\delta^{\circ}})\leq 0\mbox{ by (\ref{eq:controlconstraintsinterior}).}\]
We also have $x_{-\delta^{\circ}}\in X^{*,\epsilon}$ because $f(x_{-\delta^{\circ}})-f^*\leq \epsilon$ by assumption.

\paragraph{Second step.} Assume (\ref{eq:controlconstraintsinterior}) and also that there exists a point $z\in \widehat{X}^{*,\epsilon_0}\backslash X^{*,\epsilon}_\delta$. Since $X^{*,\epsilon}_\delta$ is closed and convex and $x_{-\delta^{\circ}}\in X^{*,\epsilon}\subset X^{*,\epsilon}_\delta$ the intersection of the line segment $[x_{-\delta^{\circ}},z]$ with $X^{*,\epsilon}_\delta$ is also closed and convex. That is, 
\[[x_{-\delta^{\circ}},z]\cap X^{*,\epsilon}_\delta = [x_{-\delta^{\circ}},x]\mbox{ with }x\in X^{*,\epsilon}_\delta.\] 
In fact we have $x\in\widehat{X}^{*,\epsilon_0+\widehat{F}(x_{-\delta^{\circ}})-\widehat{F}^*}\cap X^{*,\epsilon}_\delta$ as well. To see this, note that both $x_{-\delta^{\circ}}$ and $z$ belong to $\widehat{X}$, and this set is convex under our assumptions on $Y$ and the $f_i$, so $x\in \widehat{X}$. In addition, convexity of $\widehat{F}$ implies:
\begin{eqnarray*}\widehat{F}(x)&\leq& \max\{\widehat{F}(z),\widehat{F}(x_{-\delta^{\circ}})\}\\  
(\mbox{use that }z\in \widehat{X}^{*,\epsilon_0})&\leq &\max\{\widehat{F}^* + \epsilon_0,\widehat{F}(x_{-\delta^{\circ}})\}\\
(\mbox{note that }x_{-\delta^{\circ}}\in \widehat{X}\Rightarrow \widehat{F}(x_{-\delta^{\circ}})\geq \widehat{F}^*) &\leq &  \widehat{F}(x_{-\delta^{\circ}}) + \epsilon_0.\end{eqnarray*}

Note that $x\neq z$ and any point $x'\in (x,z]$ cannot lie in $X^{*,\epsilon}_\delta$. It follows that one of the restrictions defining $X^{*,\epsilon}_\delta$ is active at $x$. That is, one of the following properties holds:
\begin{eqnarray}\label{eq:case1xbadpoint} f(x) &=& f^*+\epsilon\mbox{ (that is, $x\in X^{*,=\epsilon}_\delta\cap \widehat{X}^{*,\widehat{F}(x_{-\delta^{\circ}})-\widehat{F}^*+\epsilon_0}$); or }\\
 \label{eq:case2xbadpoint}
\exists i\in \mathcal{I}\,:\, f_i(x) &=& \delta\mbox{ (that is, $x\in X^{*,\epsilon}_{\delta,{\rm act}(i)}\cap  \widehat{X}^{*,\widehat{F}(x_{-\delta^{\circ}})-\widehat{F}^*+\epsilon_0}$)}.\end{eqnarray}
If (\ref{eq:case1xbadpoint}) holds, then 
\[f(x) - f(x_{-\delta^{\circ}}) = \epsilon - (f(x_{-\delta^{\circ}}) -f^*)\mbox{ and }\widehat{F}(x) - \widehat{F}(x_{-\delta^\circ})\leq \epsilon_0.\]
Therefore, if (\ref{eq:case1xbadpoint}) holds, we obtain
\[\widehat{\Delta}_0(x;x_{-\delta^{\circ}}) = \widehat{\Delta}_0(x) - \widehat{\Delta}_0(x_{-\delta^\circ}) \leq \epsilon_0 -(\epsilon - (f(x_{-\delta^{\circ}}) -f^*)),\]
which means that  (\ref{eq:controlobjective}) does {\em not} hold.

Now assume (\ref{eq:case2xbadpoint}) holds. Fix an $i\in\mathcal{I}$ with $f_i(x)=\delta$. Notice that $x\in X^{*,\epsilon}_{\delta,{\rm act}(i)}$. Since $x\in\widehat{X}$ and $\widehat{F}_i(x)\leq 0$, we deduce that $\widehat{\Delta}_i(x)=\widehat{F}_i(x) - f_i(x)\leq -\delta$, which means that (\ref{eq:controlconstraintsboundary}) cannot hold. 

\paragraph{Third step.} We now assume that the three conditions in the Theorem hold. As shown above, this implies item 1 of the theorem. For simplicity, we prove item 2 assuming that some $\widehat{x}\in \widehat{X}$ achieves the minimum of $\widehat{F}$: $\widehat{F}(\widehat{x})=\widehat{F}^*$. By item 1, $\widehat{x}\in X^{*,\epsilon}_{\delta}$, so $f(\widehat{x})\leq f(x^*)+\epsilon$ and $f(x^*)-f(\widehat x)\le \gap(\delta)$. Therefore:
\begin{eqnarray*}|\widehat{F}^*-f^*|& =  & |\widehat{F}(\widehat{x}) - f(x^*)| \\ &\leq & |\widehat{F}(\widehat{x})-f(\widehat{x})| + |f(\widehat{x})- f(x^*)| \\ &\leq & \sup_{x\in X^{*,\epsilon}_{\delta}}|\widehat{\Delta}_0(x)| +  \max\{\epsilon,\gap(\delta)\} \\ 
&\leq &\sup_{x\in X^{*,\epsilon}_{\delta}}|\widehat{\Delta}_0(x;x^*)|+|\widehat{\Delta}_0(x^*)|+\max\{\epsilon,\gap(\delta)\}.
\end{eqnarray*}\qed
\end{proof}

\section{The effect of small changes in constraints on the feasible set}\label{sec:perturbation}

Our ideal optimization problem (\ref{problem:ideal}) naturally involves the feasible set $X$ and the sublevel sets $X^{*,\stheta_*}$. However, it transpires from the previous section that we will need to consider the perturbed sets $X_\delta$ and $X_{\delta}^{*,\epsilon}$, where constraints are violated by a small amount. The goal of this section is to show how one can bound the geometry and complexity of the perturbed sets in terms of the corresponding sets for the ideal problem. For this, we will make use of the geometrical assumptions from \S \ref{sub:geometry}. In what follows, $\|\cdot\|$ is a norm over $\re^d$ and $\dist$ is the corresponding set-to-point distance.

\subsection{Small constraint violations under metric regularity conditions}

The first result applies to general problems.

\begin{lemma}\label{lem:constraintviolationMRF} Make Assumption \ref{assump:MRF}. Let $\mathbb{B}$ denote the unit ball of $\re^d$ under its norm $\|\cdot\|$. Then $X_{\delta}\subset X+\mathfrak{c}\delta\,\mathbb{B}$.\end{lemma}
\begin{proof}This follows trivially from the Assumption, combined with the fact that $X_\delta\subset Y$ and the fact that $f_i(x)\leq \delta$ for all $x\in X_\delta.$\qed
\end{proof} 

\subsection{Small constraint violations under convexity}\label{sec:gamma2perturbed}

The next lemma is a key contribution of this paper. It shows that, under Assumption \ref{assump:LSCQ}, one can give a tight control of the relevant complexity parameters of $X^{*,\stheta}_{\delta}$ in terms of $X^{*,\stheta}$, for suitably small $\delta$ and $\stheta$. Recall that $x_\delta^*$ minimizes $f$ over $X_{\delta}^*$.

\begin{lemma}\label{lem:robinson2} Make Assumption \ref{assump:LSCQ}. Then:
\begin{enumerate}
\item For all $x\in Y$ with $f(x)\leq f^*+\stheta_*$ and all 
$\delta^{\circ}\in (-\eta_*,\eta_*]$, 
\[\dist(x,X^{*,\stheta_*}_{\delta^{\circ}})\leq \frac{\diam(X^{*,\stheta_*}_{\delta^{\circ}})}{\eta_*+\delta^{\circ}}\,\max_{i\in \mathcal{I}}(f_i(x)-\delta^{\circ})_+.\]
\item For all $\delta\in [0,\eta_*]$ and all $\stheta\geq {\rm gap}(-\delta)$, 
\[X^{*,\stheta}_\delta\leq 2X^{*,\stheta} - x^*_{-\delta}.\]
\item For $\delta$ and $\stheta$ as in item 2,\begin{eqnarray*}\gamma^{(\alpha)}_2(X^{*,\stheta}_\delta)&\leq &2^{\alpha}\,\gamma_2^{(\alpha)}(X^{*,\stheta})\\
\diam(X^{*,\stheta}_\delta)&\leq & 2\diam(X^{*,\stheta}).\end{eqnarray*}
\end{enumerate}\end{lemma}

The Lemma deserves some comments. Item 1 is a translation of a result of Robinson \cite{robinson1975} to our setting. Item 2 seems to be new: it states that $X^{*,\stheta}_\delta$ is contained in an homothetic copy of $X^{*,\stheta}$. This is important because, in principle, all we know from metric regularity is that $X^{*,\stheta}_\delta$ is ``close" to $X^{*,\stheta}$, meaning that $X^{*,\stheta}_{\delta}\subset X^{*,\stheta} + \mathfrak{c}\delta \mathbb{B}$ for the unit ball $\mathbb{B}$ and some constant $\mathfrak{c}>0$. By contrast, item 2 means that the actual shape of $X^{*,\stheta}_{\delta}$ is controlled by $X^{*,\stheta}$. As a result, we obtain item 3, which says that the size and complexity of $X^{*,\stheta}_{\delta}$ are controlled by the intrinsic geometry of $X^{*,\stheta}$. By contrast, one can show
\[\gamma_2^{(\alpha)}(X^{*,\stheta} + \mathfrak{c}\delta \mathbb{B})\approx \gamma_2^{(\alpha)}(X^{*,\stheta}) + c\,\delta^{\alpha}\sqrt{d}\]
for some $c>0$ depending only on $\mathfrak{c}$. In other words, metric regularity alone cannot give intrinsic bounds on the complexity of $X^{*,\stheta}_\delta$. 

We now prove the Lemma.

\begin{proof}[of Lemma \ref{lem:robinson2}]

We will need the following geometrical fact that essentially comes from Robinson's paper \cite{robinson1975}. 

\begin{claim}Take $\delta_2\in(0,\eta_*]$ and 
$\delta^{\circ}\in(-\delta_2,\delta_2]$ and $\stheta\geq f^{*}_{-\delta_2} - f^*$. Consider $x\in Y$ with $f(x)\leq f^*+\stheta$ and take $r\geq \max_{i\in\mathcal{I}}(f_i(x)-\delta^{\circ})_+\geq 0$. Let $x_{-\delta_2}^*\in X_{-\delta_2}$ be a minimizer of $f$ over that set (which exists under our assumptions of convexity and $X_{-\eta_*}\neq \emptyset$) and take \[\lambda:=\frac{r}{\delta_2 + \delta^{\circ} +r}\in [0,1).\] Then $x^{(\lambda)} := (1-\lambda)x + \lambda x^*_{-\delta_2}\in X^{*,\stheta}_{\delta^{\circ}}.$ \end{claim}

Indeed, it is obvious that $x^{(\lambda)}\in Y$ because this set is convex. We also have that $f(x^{(\lambda)})\leq f^*+\stheta$ because $f$ is convex and both $x$, $x^*_{-\delta_2}$ satisfy this inequality. Finally, for each $i\in\mathcal{I}$, 
\[f_i(x^{(\lambda)})-\delta^{\circ} \leq (1-\lambda)(f_i(x)-\delta^{\circ}) + \lambda (f_i(x_{-\delta_2}^*)-\delta^{\circ})\leq (1-\lambda)r -\lambda (\delta_2+\delta^{\circ}) =0.\]So $x^{(\lambda)}\in X^{*,\stheta}_{\delta^{\circ}}$. 

We now use this Claim to obtain parts 1 and 2 of the Lemma. We will then obtain part 3 from part 2.

\paragraph{Proof of Lemma \ref{lem:robinson2}, part 1.} We apply the claim to $x$ as in item $1$ with $\delta_2=\eta_*$, $r:=\max_{i\in\mathcal{I}}(f_i(x)-\delta^{\circ})_+$ and $\stheta=\stheta_*$. In that case, we see that:
\[x^{(\lambda)}\in X^{*,\stheta}_{\delta^{\circ}}\Rightarrow \dist(x,X_{\delta^{\circ}}^{*,\stheta})\leq \|x - x^{(\lambda)}\|.\]
Since
\[x - x^{(\lambda)} = \lambda (x - x^{*}_{-\eta_*}) \mbox{ and } x^{(\lambda)}- x^{*}_{-\eta_*} =  (1-\lambda)(x-x^*_{-\eta_*}),\]
\[\|x - x^{(\lambda)}\| = \frac{\lambda}{1-\lambda} \|x^{(\lambda)}-x^*_{-\eta_*}\|\leq \frac{r}{\delta_2+\delta^{\circ}}\,\diam(X^{*,\stheta}_{\delta^{\circ}})\] because both $x^*_{-\eta_*}$ and $x^{(\lambda)}$ belong to $X^{*,\stheta}_{\delta^{\circ}}.$ Noting that $\delta_2=\eta_*$, gives the result. 

\paragraph{Proof of Lemma \ref{lem:robinson2}, part 2.} For $\delta=0$ the claim is trivial. Suppose $\delta\neq0$. Take an arbitrary $x \in X^{*,\stheta}_{\delta}$; in particular, $f_0(x)\leq f^*+\stheta$ and $f_i(x)\leq \delta$ for all $i\in\mathcal{I}$. Apply the claim with $\delta^{\circ}=0$, $r=\delta$ and $\delta_2=\delta$. In this case, $\lambda=1/2$ and therefore,
\[\frac{x+x^*_{-\delta}}{2} = x^{(\lambda)}\in X^{*,\stheta}\Rightarrow x = 2x^{(\lambda)} - x^*_{-\delta}\in 2X^{*,\stheta} - x^{*}_{-\delta}.\]

\paragraph{Proof of Lemma \ref{lem:robinson2}, part 3.} We combine the previous item with the following simple facts. The first is that $\gamma_2^{(\alpha)}$ and $\diam$ are invariant under translations. Moreover, $\gamma_2^{(\alpha)}(\lambda S) = \lambda^{\alpha}\gamma_2^{(\alpha)}(S)$ and $\diam(\lambda S) =\lambda\,\diam(S)$ for all $S\subset\re^d$ and $\lambda\geq 0$.
\qed
\end{proof}

\section{Proofs of main results}\label{sec:proofs}

We combine here the tools from the previous three sections to prove Theorems \ref{thm:generalrandomset} (in \S \ref{proof:generalrandomset}) and \ref{thm:convexrandomset} (in \S \ref{proof:convexrandomset}). 

\subsection{General sets and functions}\label{proof:generalrandomset}

\begin{proof}[of Theorem \ref{thm:generalrandomset}] The strategy of the proof is as follows. We will use Lemma \ref{lem:generaldeviations} to show that the event ${\rm Good}_{\rm Thm. \ref{thm:generalrandomset}}(t,\epsilon_0)$ contains the intersection of events $E_1$ and $E_2$ below. We then lower bound 
$\prob(E_1)$ and $\prob(E_2)$ to finish the proof.

The two events are defined as follows. 
\begin{eqnarray}E_1 &:=& \bigcap_{i\in\mathcal{I}_0}\left\{|\widehat{\Delta}_i(x^*)|\leq \sigma_*\sqrt{\frac{6(1+\log (2|\mathcal{I}|+2) + t)}{N}}\right\};\\
E_2&:=& \bigcap_{i\in\mathcal{I}_0}\left\{\sup_{x\in Y}|\widehat{\Delta}_i(x;x^*)|\leq \widehat{r}_N(t)\right\}.\end{eqnarray}

\paragraph{First part: containement.} Recalling the definitions of $\widehat{\Delta}_i$ in Section \ref{sec:deviationlocalization}, we see that $\widehat{\Delta}_i(y) = \widehat{\Delta}_i(y;x^*) + \widehat{\Delta}_i(x^*)$ for all $y\in X$, so:
\[\sup_{y\in Y,i\in\mathcal{I}_0}|\widehat{\Delta}_i(y)|\leq \sup_{y\in Y,i\in\mathcal{I}_0}|\widehat{\Delta}_i(y;x^*)| + |\widehat{\Delta}_i(x^*)|.\]
In particular, 
\[\mbox{if $E_1\cap E_2$ holds, }\sup_{y\in Y,i\in \mathcal{I}_0}|\widehat{\Delta}_i(y)|\leq \sigma_*\sqrt{\frac{6(1+\log (2|\mathcal{I}|+2) + t }{N}} + \widehat{r}_N(t).\]
So the assumptions of Lemma \ref{lem:generaldeviations} are satisfied with 
\[\delta:=\sigma_*\sqrt{\frac{6(1+\log (2|\mathcal{I}|+2) + t }{N}} + \widehat{r}_N(t)\mbox{ and }\epsilon:= \widehat{r}_N(t).\]
Applying the Lemma and inspecting the definitions shows that ${\rm Good}_{\rm Thm. \ref{thm:generalrandomset}}(t,\epsilon_0)$ holds. Indeed, We deduce that $E_1\cap E_2\subset {\rm Good}_{\rm Thm. \ref{thm:generalrandomset}}(t,\epsilon_0)$ holds. 

\paragraph{Second part: probability bounds} To finish, we must prove that $\prob(E_1\cap E_2)\geq 1 - e^{-t} - (2|\mathcal{I}|+1)\rho$. Note that:
\begin{eqnarray}\label{eq:line1probboundsgeneralrandom}1 - \prob(E_1\cap E_2) &\leq & \sum_{i\in\mathcal{I}_0}\prob\left\{|\widehat{\Delta}_i(x^*)|>\sigma_*\sqrt{\frac{6(1+\log (2|\mathcal{I}|+2) + t }{N}}\right\} \\ \label{eq:line2probboundsgeneralrandom}& & + \sum_{i\in\mathcal{I}_0}\prob\left\{\sup_{x\in Y}|\widehat{\Delta}_i(x;x^*)|> \widehat{r}_N(t)\right\} \end{eqnarray}

Therefore, it suffices to bound each term in (\ref{eq:line1probboundsgeneralrandom}) and (\ref{eq:line2probboundsgeneralrandom}) separately. For the terms in (\ref{eq:line1probboundsgeneralrandom}), we apply Lemma \ref{lem:panchenko} with $\nu^2=\sigma_*^2$ and $Z_k:= \pm F_i(x^*,\xi_k)$, so that:
\[\frac{1}{N}\sum_{k=1}^N(Z_k - \esp[Z_1]) = \pm (\widehat{\espn} - \espn)F_i(x^*,\cdot).\]
Because of Assumption \ref{assump:goodFi}, we know that:
\[\prob\left\{\frac{1}{N}\sum_{k=1}^N(Z_k - \esp[Z_1])^2\geq 2\sigma_*^2\right\}\leq \rho.\]
Therefore Lemma \ref{lem:panchenko} gives:
\[\prob\left\{|(\widehat{\espn} - \espn)F_i(x^*,\cdot)|>\sigma_*\sqrt{\frac{6(1+\log (2|\mathcal{I}|+2) + t)}{N}}\right\} \leq   \frac{e^{-t}}{2(|\mathcal{I}|+2)} + \rho.\]
To bound the terms in (\ref{eq:line2probboundsgeneralrandom}) we fix an $i\in\mathcal{I}_0$ and apply our concentration result, Theorem \ref{thm:concentration}. In the language of that theorem, we have
\[\widehat{\Delta}_i(x;x^*) = (\widehat{\espn} - \espn)\,(G(x,\cdot) - G(x^*,\cdot))\mbox{ for }G:=F_i.\]
With these choices, 
\[\Delta G = \sup_{x\in Y}|\widehat{\Delta}_i(x;x^*)|.\]
Assumption \ref{assump:goodFi} guarantees that:
\[\prob\left\{\widehat{\espn}\mathsf{L}_i(\cdot)>2\sigma^2\right\}\leq \rho.\]
So Theorem \ref{thm:concentration} is applicable with $\nu^2=\sigma^2$, $\mathsf{L} = \mathsf{L}_i$. Checking the formula for $\widehat{r}_N(t)$, we may now use Theorem \ref{thm:concentration} to deduce:
\[\prob\left\{\sup_{x\in Y}|\widehat{\Delta}_i(x;x^*)|>\widehat{r}_N(t)\right\}\leq \frac{e^{-t}}{2(|\mathcal{I}|+2)} + \rho.\]
We have now bounded all the terms in the sums (\ref{eq:line1probboundsgeneralrandom}) and (\ref{eq:line2probboundsgeneralrandom}). Plugging the bounds back into these equations give the desired lower bound on $\prob(E_1\cap E_2)$.
\qed
\end{proof}

\subsection{Convex sets and functions}\label{proof:convexrandomset}

\begin{proof}[of Theorem \ref{thm:convexrandomset}]For convenience, we only consider the case where $\mathcal{I}\neq \emptyset$, as the other case is simpler.

Our general proof strategy is similar to the one of Theorem \ref{thm:generalrandomset}. In the first step of the proof, we define {\em decreasing} sequences of events $E_{1,k},E_{2,k}$ and argue that ${\rm Good}_{\rm Thm. \ref{thm:convexrandomset}}(t,\epsilon_0)$ contains $\cap_{k}(E_{1,k}\cap E_{2,k})$. We then bound the probability of the good event via bounds on $\prob(E_{1,k})$ and $\prob(E_{2,k})$. 

Let us first define the events. Looking at the definition of $\check{r}(t;\epsilon_0)$, we see that one can find a {\em decreasing} sequence $\{\epsilon_k\}_{k=1}^{+\infty}$.
\begin{equation*}\forall k\geq 1 \,:\,\epsilon_k \in R_{N,\eta_*}(t;\epsilon_0)\mbox{ and moreover }\epsilon_k\searrow \check{r}(t;\epsilon_0).\end{equation*}
For each $k$, we have:
\[2\check{w}(t;\epsilon_k) + {\rm gap}(-\check{\delta}(t;\epsilon_k)) +\epsilon_0 <\epsilon_k.\]
Now, $\check{w}(t;\epsilon_k)=\widehat{w}_N(t;\check{\delta}(t;\epsilon_k);\epsilon_k)$ where $\check{\delta}(t;\epsilon_k) = \inf S_{N,\eta_*}(t;\epsilon_k)$. Given our assumptions, (\ref{eq:holdercont}) and Lemma \ref{lem:robinson2} above, it is easy to check that 
\[\delta \mapsto \widehat{w}_N(t;\delta;\epsilon_k)\mbox{ and }\delta \mapsto {\rm gap}(-\delta)\]
are continuous nonincreasing functions of $\delta\in [0,\eta_*]$. Moreover, the sets $S_{N,\eta_*}(t;\epsilon_k)$ increase with $k$, so $\check{\delta}(t;\epsilon_k)$ decreases with $k$. Therefore, one can find a {\em decreasing} sequence $\{\delta_k\}_{k=1}^{+\infty}$ such that:
\begin{eqnarray}\nonumber \forall k\geq 1& : & \delta_k\in S_{N,\eta_*}(t;\epsilon_k),\\ \nonumber & & 2\widehat{w}_N(t;\delta_k;\epsilon_k) + {\rm gap}(-\delta_k) +\epsilon_0 <\epsilon_k,\\ \label{eq:defdeltak} & &  \mbox{ and }\delta_k<\check{\delta}(t;\epsilon_k)+k^{-1}.\end{eqnarray}
It follows in particular, that 
\[\lim\delta_k = \lim_{k}\check{\delta}(t;\epsilon_k) = \lim_{\epsilon\searrow \check{r}_N(t;\epsilon_0)}\check{\delta}(t;\epsilon)=\check{\delta}(t).\]

The events we define are:
\begin{eqnarray}\label{eq:defE1convex}E_{1} &:=& \bigcap_{i\in\mathcal{I}_0}\left\{|\widehat{\Delta}_i(x^*)|\leq \sigma_*\sqrt{\frac{6(1+\log (2|\mathcal{I}|+2) + t)}{N}}\right\};\\ \label{eq:defE2convex}
E_{2,k}&:=& \bigcap_{i\in\mathcal{I}_0}\left\{\sup_{x\in X_{\delta_k}^{*,\epsilon_k+{\rm gap}(-\delta_k)}}|\widehat{\Delta}_i(x;x^*)|\leq \widehat{w}_N(t;\delta_k;\epsilon_k)\right\}.\end{eqnarray}

The fact that $\{\epsilon_k\}_{k}$ and $\{\delta_k\}$ are both decreasing implies that the events $E_{2,k}$ are decreasing.

\paragraph{First part: containment.} We will argue that $\cap_{k}(E_{1}\cap E_{2,k})\subset {\rm Good}_{\rm Thm. \ref{thm:convexrandomset}}(t,\epsilon_0)$. To show this, we assume that the event $\cap_{k}(E_{1}\cap E_{2,k})$ holds, and deduce that ${\rm Good}_{\rm Thm. \ref{thm:convexrandomset}}(t,\epsilon_0)$ must hold as well. 

Fix an index $k$. We may assume  that there exists $x^{*}_{-\delta_k}$ minimizing $f$ over $X_{-\delta_k}$ and:
\[f(x^*_{-\delta_k}) - f^* = {\rm gap}(-\delta_k).\]
In particular, $x^*_{-\delta_k}\in X^{*,\epsilon_k}_{\delta_k}$. Because $E_{1}\cap E_{2,k}$ holds, we have that for all $i\in\mathcal{I}_0$ and $x\in X^{*,\epsilon_k+{\rm gap}(-\delta_k)}_{\delta_k}$:
\begin{eqnarray*}|\widehat{\Delta}_i(x)| &\leq & |\widehat{\Delta}_i(x;x^*)| + |\widehat{\Delta}_i(x^*)|\\ &\leq & \widehat{w}_N(t;\delta_k;\epsilon_k) +\sigma_*\sqrt{\frac{6(1+\log (2|\mathcal{I}|+2) + t)}{N}}\\ (\mbox{use }\delta_k\in S_{N,\eta_*}(t;\epsilon_k)) &<& \delta_k.\end{eqnarray*}
In particular, $\widehat{\Delta}_i(x^*_{-\delta_k})<\delta_k$. For the same $x$'s, we also have:
\begin{eqnarray*}|\widehat{\Delta}_0(x;x_{-\delta_k}^*)|&\leq & |\widehat{\Delta}_0(x^*_{-\delta_k};x^*)| +  |\widehat{\Delta}_0(x;x^*)|\\ &\leq & 2\widehat{w}_N(t;\delta_k;\epsilon_k) \\ \mbox{(use (\ref{eq:defdeltak}))}&<&  \epsilon_k - \epsilon_0 - {\rm gap}(-\delta_k).\end{eqnarray*}

We now apply Proposition \ref{lem:convexlocalized} with $\delta^{\circ} = \delta=\delta_k$, $x_{-\delta^{\circ}}:=x_{-\delta_k}^*$ and $\epsilon:=\epsilon_k$. The above calculations imply that the three conditions of such lemma, given by  (\ref{eq:controlconstraintsinterior}), (\ref{eq:controlconstraintsboundary}) and (\ref{eq:controlobjective}), are satisfied. We conclude that
\begin{equation}\label{eq:forallk1}\widehat{X}^{*,\epsilon_0}\subset X_{\delta_k}^{*,\epsilon_k}\end{equation} 
and (by the same estimates)
\begin{eqnarray}\nonumber |\widehat{F}^* - f^*|&\leq & |\widehat{\Delta}_0(x^*)|  + \sup_{x\in X^{*,\epsilon_k}_{\delta_k}}|\widehat{\Delta}_0(x;x^*)| + \max\{\epsilon_k,\gap(\delta_k)\} \\ \nonumber (\mbox{$E_1$ holds}) &\leq &  \sigma_*\sqrt{\frac{6(1+\log (2|\mathcal{I}|+2) + t)}{N}} + \sup_{x\in X^{*,\epsilon_k}_{\delta_k}}|\widehat{\Delta}_0(x;x^*)| + \max\{\epsilon_k,\gap(\delta_k)\}\\ \nonumber \mbox{($E_{2,k}$ occurs)} &\leq & \sigma_*\sqrt{\frac{6(1+\log (2|\mathcal{I}|+2)  + t)}{N}} + \widehat{w}_N(t;\delta_k,\epsilon_k) + \max\{\epsilon_k,\gap(\delta_k)\} \\ \label{eq:forallk2} \mbox{(use (\ref{eq:defdeltak}))} &\leq &  \sigma_*\sqrt{\frac{6(1+\log (2|\mathcal{I}|+2)  + t)}{N}} + \frac{\epsilon_k}{2}+\max\{\epsilon_k,\gap(\delta_k)\}.
\end{eqnarray}

Both (\ref{eq:forallk1}) and (\ref{eq:forallk2}) hold for all $k$. Letting $k\to +\infty$ and recalling that  $\delta_k\searrow \check{\delta}(t)$ and $\epsilon_k\searrow \check{r}(t;\epsilon_0)$, we obtain:
\begin{eqnarray*}\label{eq:limit1}\widehat{X}^{*,\epsilon_0}&\subset & X_{\check{\delta}(t)}^{*,\check{r}(t;\epsilon_0)};\\  |\widehat{F}^* - f^*| &\leq & \sigma_*\sqrt{\frac{6(1+\log (2|\mathcal{I}|+2)  + t)}{N}} + \frac{\check{r}(t;\epsilon_0)}{2}
+\max\{\check{r}(t;\epsilon_0),\gap(\check\delta(t))\}.
\end{eqnarray*}

Going back to the statement of Theorem \ref{thm:convexrandomset} (page \pageref{thm:convexrandomset}), we see that the two properties above correspond to {\bf (a)} and {\bf (b)} in the definition of ${\rm Good}_{\rm Thm. \ref{thm:convexrandomset}}(t,\epsilon_0)$. The remaining property {\bf (c)} that defines that event also holds due to Lemma \ref{lem:robinson2}. Therefore, by assuming that $E_1\cap E_{2,k}$ occurs for all $k$, we have deduced that ${\rm Good}_{\rm Thm. \ref{thm:convexrandomset}}(t,\epsilon_0)$ also holds.

\paragraph{Second step: probability bounds} Recall that the events $E_{2,k}$ are decreasing. By the first step, 
\[\prob({\rm Good}_{\rm Thm. \ref{thm:convexrandomset}}(t,\epsilon_0))\geq \prob\left(E_1\cap \bigcap_{k=1}^{+\infty} E_{2,k}\right) = \lim_{k\to +\infty}\prob(E_1\cap E_{2,k}).\]
Therefore, all that remains to show is that:
\[{\bf Goal: }\forall k\geq 1\,:\, 1 - \prob(E_1\cap  E_{2,k})\leq e^{-t} + 2(|\mathcal{I}|+1)\rho.\]
From this point on, the proof resembles the second step in the proof of Theorem \ref{thm:generalrandomset}, and we will be a bit briefer. Following (\ref{eq:line1probboundsgeneralrandom}) and (\ref{eq:line2probboundsgeneralrandom}), but with the definition of $E_{2,k}$ in (\ref{eq:defE2convex}), we obtain 
\begin{eqnarray*}  \prob(E^c_1\cup E^c_{2,k}) &\leq & \sum_{i\in\mathcal{I}_0}\prob\left\{|\widehat{\Delta}_i(x^*)|>\sigma_*\sqrt{\frac{6(1+\log (2|\mathcal{I}|+2) + t }{N}}\right\} \\ & & + \sum_{i\in\mathcal{I}_0} \prob\left\{\sup_{x\in X_{\delta_k}^{*,\epsilon_k+{\rm gap}(-\delta_k)}}|\widehat{\Delta}_i(x;x^*)|>\widehat{w}_N(t;\delta_k;\epsilon_k)\right\}. \end{eqnarray*}
As in the proof of Theorem \ref{thm:generalrandomset}, Lemma \ref{lem:panchenko} gives:
\[\forall i\in\mathcal{I}_0\,:\, \prob\left\{|\widehat{\Delta}_i(x^*)|>\sigma_*\sqrt{\frac{6(1+\log (2|\mathcal{I}|+2) + t }{N}}\right\}\leq \frac{e^{-t}}{2(|\mathcal{I}|+1)} + \rho.\]
On the other hand, the bound
\[\forall i\in\mathcal{I}_0\,:\,\prob\left\{\sup_{x\in X_{\delta_k}^{*,\epsilon_k+{\rm gap}(-\delta_k)}}|\widehat{\Delta}_i(x;x^*)|>\widehat{w}_N(t;\delta_k;\epsilon_k)\right\}\leq \frac{e^{-t}}{2(|\mathcal{I}|+1)} + \rho\]
follows from applying Theorem \ref{thm:concentration} as in the proof of Theorem \ref{thm:generalrandomset}, noting that this time we have Assumption \ref{assump:goodFi} over $Z=X^{*,\stheta_*}_{\eta_*}\supset X_{\delta_k}^{*,\epsilon_k+{\rm gap}(-\delta_k)}$, and also that 
\[\gamma^{(\alpha)}_2(X_{\delta_k}^{*,\epsilon_k+{\rm gap}(-\delta_k)})\leq 2\gamma^{(\alpha)}_2(X^{*,\epsilon_k+{\rm gap}(-\delta_k)})\]
and \[{\diam}(X_{\delta_k}^{*,\epsilon_k+{\rm gap}(-\delta_k)})\leq 2{\diam}(X^{*,\epsilon_k+{\rm gap}(-\delta_k)})\]
by Lemma \ref{lem:robinson2}.
\qed
\end{proof}

\section*{Appendix}

\begin{proof}[of Lemma \ref{lem:panchenko}] The second statement in the Lemma is a direct consequence of the first. Therefore, we will only prove the first statement. 

Assume that $Z'_1,\dots,Z'_n$ are independent copies of the $Z_1,\dots,Z_n$. Also let $Z=(Z_1,\dots,Z_n)^T$. What we want to prove is that, for any $t\geq 0$, 
\[{\bf Want:}\;\prob\left\{\esp\left[\sum_{k=1}^N (Z_k - Z'_k)\mid Z\right]\geq \sqrt{2(1+t)\,\sum_{k=1}^N\esp[(Z_k-Z'_k)^2\mid Z]} \right\}
\leq e^{-t}.\]
By \cite[Corollary 1]{panchenko}, it suffices to prove that, for any $t\geq 0$,
\[{\bf Sufficient:}\;\prob\left\{\sum_{k=1}^N (Z_k - Z'_k)\geq \sqrt{2t\,\sum_{k=1}^N(Z_k-Z'_k)^2} \right\}\leq e^{-t}.\]
We will prove that the above inequality holds almost surely conditionally on values $|Z_k-Z'_k|=a_k$, $1\leq k\leq N$. Notice that, conditionally on these values,
\[Z_k - Z'_k = u_k\,a_k\]
where the $u_k$ are i.i.d. unbiased random signs. So what we must show is that:
\[\forall t\geq 0\,:\, \prob\left\{\sum_{k=1}^N u_ia_i\geq\sqrt{2t \sum_{k=1}^Na_k^2}\right\}\leq e^{-t},\]
for any choice of $a_k$, $1\leq k\leq N$. This follows easily from the standard inequalities:
\[\forall \theta>0\,:\,\esp[e^{\theta\sum_{k=1}^N u_ia_i}] = \prod_{k=1}^N\cosh(\theta a_k)\leq e^{\frac{\theta^2\sum_{k=1}^Na_k^2}{2}},\]
and Bernstein's trick:
\[\prob\left\{\sum_{k=1}^N u_ia_i\geq\sqrt{2t \sum_{k=1}^Na_k^2}\right\}\leq \inf_{\theta>0}\esp[e^{\theta\sum_{k=1}^N u_ia_i}]e^{-\theta\sqrt{2t \sum_{k=1}^Na_k^2}}\leq e^{-t}.\]
\qed
\end{proof}

\begin{proof}[of Proposition \ref{prop:goodisgreath}] We will need the following Lemma. 

\begin{lemma}\label{lem:BDG}There exists a constant $\cbdg$ such that, for all $p\geq 2$ and all i.i.d. random variables $Z_1,\dots,Z_N\in L^p$ with $\esp[Z_i]=0$,
\[\Lpnorm{\frac{Z_1 + \dots + Z_N}{N}}\leq  \cbdg\sqrt{\frac{p}{N}}\,\Lpnorm{Z_1},\]\end{lemma}
\begin{proof}[of the Lemma] By the Burkholder-Davis-Gundy inequality and the subaditivity of the $L^{p/2}$ norm:
\[\Lpnorm{Z_1 + \dots + Z_N}\leq \cbdg\sqrt{p}\,\Lunorm{Z^2 _1+\dots+Z^2_N}{p/2}^{1/2}\leq  \cbdg\sqrt{p\sum_{i=1}^N\Lunorm{Z^2_i}{p/2}}\] 
and the proof finishes when we note $\Lunorm{Z^2_i}{p/2} =\Lpnorm{Z_1}^2$ for each index $i$.\qed
\end{proof}

Now note that the random variables \[H_k:=\frac{h(\xi_k) - \espn h(\cdot)}{\sigma^2}\,\,(1\leq k\leq N)\]
are i.i.d. and satisfy $\esp[H_k]=0$, $\Lpnorm{H_k}\leq \kappa_p$. Markov's inequality implies:
\[\prob\left\{\widehat{\espn}h(\cdot)> 2\sigma^2\right\}\le \prob\left\{\frac{1}{N}\sum_{k=1}^NH_k>1\right\}\leq \Lpnorm{\frac{1}{N}\sum_{k=1}^NH_k}^p.\]
Now use Lemma \ref{lem:BDG} to bound the RHS.\qed \end{proof}

\end{document}